\documentclass[12pt]{amsart}
\usepackage{amsfonts,amsmath,amssymb}
\newtheorem{LL}{LL}[section]

\newtheorem{Prop}[LL]{Proposition}
\newtheorem{Cor}[LL]{Corollary}

\newtheorem{Theorem}[LL]{Theorem}
\newtheorem{Lemma}[LL]{Lemma}
\theoremstyle{definition}

\newtheorem{Ex}[LL]{Example}

\newcommand{\pf}{\medskip\noindent{\sc Proof: }}
\newcommand{\epf}{$\Box$}

\newcommand{\PH}{\!\!\!\!\!\!\!\!\!\!\!\!\!\! \phantom{\begin{array}{r} a /a \end{array}}}

\newcommand\Alg{\operatorname{Alg}}

\newcommand\id{{\operatorname{id}}}

\newcommand\co{{\operatorname{co}}}

\newcommand\Irr{{\operatorname{Irr}}}

\newcommand\D{\mathcal{D}}

\newcommand\G{\Gamma}

\newcommand\sw[1]{{}_{(#1)}}

\title{Representations Parameterized by a Pair of Characters}
\author{David E. Radford}
\thanks{Research by the first author partially supported by NSA Grant H98230-04-1-0061. A
significant amount of work on this paper was done during his
visits to the Mathematisches Institut der
Ludwig-Maximilians-Universit\"{a}t M\"{u}nchen during June of 2003
and May of 2004 and during the visits of the second author to the
Department of Mathematics, and Statistics, and Computer Science at
the University of Illinois at Chicago during September 2003 and
March 2005. The first author expresses his gratitude for the
hospitality and support he received from the Institut and the
second expresses his gratitude for the same he received from UIC}
\address{University of Illinois at Chicago \\
Department of Mathematics, Statistics and \\ Computer Science (m/c 240) \\
801 South Morgan Street \\
Chicago, IL   60607-7045} \email{radford@uic.edu}
\author{Hans J\"{u}rgen Schneider}
\address{Mathematisches Institut \\
Ludwig-Maximilians-Universit\"{a}t M\"{u}nchen \\
Theresienstr. 39 \\
D-80333 M\"{u}nchen, Germany}
\email{Hans-Juergen.Schneider@mathematik.uni-muenchen.de}
\numberwithin{equation}{section}
\begin{document}

\maketitle

\date{}
\begin{abstract}
{ \small  \rm Let $U$ and $A$ be algebras over a field $k$. We
study algebra structures $H$ on the underlying tensor product
$U{\otimes}A$ of vector spaces which satisfy
$(u{\otimes}a)(u'{\otimes}a') = uu'{\otimes}aa'$ if $a = 1$ or $u'
= 1$. For a pair of characters $\rho \in \Alg(U, k)$ and $\chi \in
\Alg(A, k)$ we define a left $H$-module $L(\rho, \chi)$. Under
reasonable hypotheses the correspondence $(\rho, \chi) \mapsto
L(\rho, \chi)$ determines a bijection between character pairs and
the isomorphism classes of objects in a certain category
${}_H\underline{\mathcal M}$ of left $H$-modules. In many cases
the finite-dimensional objects of ${}_H\underline{\mathcal M}$ are
the finite-dimensional irreducible left $H$-modules.

In \cite{Irr} we apply the results of this paper and show that the
finite-dimensional irreducible representations of a wide class of
pointed Hopf algebras are parameterized by pairs of characters.}
\end{abstract}
\setcounter{section}{0}
\section*{Introduction}\label{SecIntro}
This paper develops the theory of a type of modules for certain
algebra structures $H$ defined on tensor products which, in many
cases, accounts for the finite-dimensional irreducible
representations of $H$. The general results are applied to the
study of irreducible modules of a certain class of pointed Hopf
algebras over a field $k$ in \cite{Irr}. This class is very basic
in light of the results of the program of Andruskiewitsch and the
second author to determine the structure of
 pointed Hopf algebras with commutative
coradicals \cite{AS3,AS2,SmallAS,ASSurvey}. The Hopf algebras of
interest to us are quotients of certain two-cocycle twists $H =
(U{\otimes}A)^\sigma$ of the tensor product of two pointed Hopf
algebras $U$ and $A$ over $k$. As a vector space $H = U{\otimes}A$
and multiplication has the property $(u{\otimes}a)(u'{\otimes}a')
= uu'{\otimes}aa'$ whenever $a = 1$ or $u' = 1$.

The natural context for us to begin our study is in the category
of algebra structures on the tensor product $U{\otimes}A$ of the
underlying vector spaces of algebras $U$ and $A$ over $k$ which
satisfy the multiplication property. For a pair of characters
$\rho \in \Alg (U, k)$ and $\chi \in \Alg (A, k)$ we construct a
left $H$-module $L(\rho, \chi)$ and a right $H$-module $R(\chi,
\rho)$. They satisfy a duality relationship with respect to a
certain $H$-balanced bilinear form $\Psi$.

We describe the modules $L(\rho, \chi)$ and $R(\chi, \rho)$
abstractly. The $H$-module $L(\rho, \chi)$ contains a codimension
one left $U$-submodule and is generated as an $H$-module by a
one-dimensional $A$-submodule. Let ${}_H\underline{\mathcal M}$ be
the full subcategory of all left $H$-modules whose objects contain
a left $U$-submodule of codimension one and contain a
one-dimensional left $A$-submodule. Under mild conditions we show
the correspondence $(\rho , \chi) \mapsto L(\rho, \chi)$ is
one-to-one. One of our main results gives natural conditions under
which this correspondence determines a bijection between the
Cartesian product $\Alg(U, k){\times}\Alg(A, k)$ and the
isomorphism classes of ${}_H\underline{\mathcal M}$. We will need
to know when finite-dimensional irreducible representations are
one-dimensional.

This paper is organized as follows. In Section \ref{SecPrelim} we
set our notation for algebras, coalgebras, Hopf algebras, and the
like. Two-cocycles $\sigma$ are reviewed and the Hopf algebra $H =
(U{\otimes}A)^\sigma$ is described. The Drinfeld double of a
finite-dimensional Hopf algebra $H$ over $k$ can be viewed as
$D(H) = (H^{* \, cop}{\otimes}H)^\sigma$ for some two-cocycle.

Let $A$ and $U$ be algebras over $k$. In Section \ref{SecUOtimesA}
we develop a theory of algebra structures $H$ on $U{\otimes}A$
which satisfy the multiplication property. We define the
$H$-modules $L(\chi, \rho)$, $R(\rho, \chi)$ and study them as
explicit constructions and also more abstractly. The $H$-modules
$L(\chi, \rho)$ are objects of ${}_H\underline{\mathcal M}$, the
category whose objects $M$ are left $H$-modules which contain a
codimension one left $U$-submodule $N$ and are generated as an
$H$-module by a one-dimensional left $A$-submodule $km$. We study
this category and its refinement ${}_H\underline{\mathcal M}'$,
whose objects are triples $(M, N, km)$, and duality relations with
their counterparts for right $H$-modules. There are natural
$H$-balanced bilinear forms which provide duality relationships
between the objects of ${}_H\underline{\mathcal M}'$ and
$\underline{\mathcal M}'_H$.

Section \ref{SecMainHM} contains the main results for the modules
$L(\chi, \rho)$ in the more abstract setting of the category
${}_H\underline{\mathcal M}$. We first consider conditions on the
characters $\rho$ and $\chi$ separately in connection with the
correspondence $(\rho, \chi) \mapsto L(\rho, \chi)$. Our main
theorem gives natural conditions under which this correspondence
determines a bijection between $\Alg(U, k){\times}\Alg(A, k)$ and
the isomorphism classes of the objects of ${}_H\underline{\mathcal
M}$. In this case each object of ${}_H\underline{\mathcal M}$ has
a unique codimension one $U$-submodule and a unique
one-dimensional left $A$-submodule. Under the hypothesis of the
theorem, when the finite-dimensional irreducible representations
of $U$ and $A$ are one-dimensional then the finite-dimensional
irreducible left $H$-modules are the same as the
finite-dimensional objects of ${}_H\underline{\mathcal M}$. We
consider the case when $U$ and $A$  are Hopf algebras in Section
\ref{SecMainBialg}. In Section \ref{SecFinDimOne} we consider
conditions under which the irreducible representations of an
algebra are one-dimensional for applications to certain classes of
pointed Hopf algebras.

Throughout $k$ is a field and all vector spaces are over $k$. For
vector spaces $U$ and $V$ we will drop the subscript $k$ from
${\rm End}_k(V)$, ${\rm Hom}_k(U, V)$, and $U{\otimes}_kV$. We
denote the identity map of $V$ by $\id_V=\id$. For a non-empty
subset $S$ of the dual space $V^*$ we let $S^\perp$ denote the
subspace of $V$ consisting of the common zeros of the functionals
in $S$. The ``twist" map $\tau_{U, V} : U{\otimes}V
\longrightarrow V{\otimes}U$ is defined by $\tau_{U,
V}(u{\otimes}v) = v{\otimes}u$ for all $u \in U$ and $v \in V$.
For $p \in U^*$ and $u \in U$ we denote the evaluation of $p$ on
$u$ by $p(u)$ or ${<}p, u{>}$. Any one of
\cite{LamRad,Mont,SweedlerBook} will serve as a Hopf algebra
reference for this paper.
\section{Preliminaries}\label{SecPrelim}
For a group $G$ we let $\widehat{G}$ denote the group of
characters of $G$ with values in $k$. $H = kG$ denotes the group
algebra of $G$ over $k$ which is a Hopf algebra arising in most
applications in this paper. We usually denote the antipode of a
Hopf algebra over $k$ by $S$.

Let $(A, m, \eta)$ be an algebra over $k$, which we shall usually
denote by $A$. Generally we represent algebraic objects defined on
a vector space by their underlying vector space. We say that $a, b
\in A$ {\em skew commute} if $ab = \omega ba$ for some non-zero
$\omega \in k$. Note that  $(A, m^{op}, \eta)$ is an algebra over
$k$, where $m^{op} = m{\circ}\tau_{A, A}$. We denote $A$ with this
algebra structure by $A^{op}$ and we denote the category of left
(respectively right) $A$-modules and module maps by ${}_A{\mathcal
M}$ (respectively ${\mathcal M}_A$). If ${\mathcal C}$ is a
category, by abuse of notation we will write $C \in {\mathcal C}$
to indicate that $C$ is an object of ${\mathcal C}$.

Let $M$ be a left $A$-module. Then $M^*$ is a right $A$-module
under the transpose action  which is given by $(m^*{\cdot}a)(m) =
m^*(a{\cdot}m)$ for all $m^* \in M^*$, $a \in A$, and $m \in M$.
Likewise if $M$ is a right $A$-module then $M^*$ is a left
$A$-module where $(a{\cdot}m^*)(m) = m^*(m{\cdot}a)$ for all $a
\in A$, $m^* \in M^*$, and $m \in M$.

Let $(C, \Delta, \epsilon)$ be a coalgebra over $k$, which we
usually denote by $C$. Generally we use a variant on the
Heyneman-Sweedler notation for the coproduct and write $\Delta (c)
= c_{(1)}{\otimes}c_{(2)}$ to denote $\Delta(c) \in C{\otimes}C$
for $c \in C$. Note that $(C, \Delta^{cop}, \epsilon)$ is a
coalgebra over $k$, where $\Delta^{cop} = \tau_{C,
C}{\circ}\Delta$. We let $C^{cop}$ denote the vector space $C$
with this coalgebra structure. Observe that $C$ is a
$C^*$-bimodule with the actions defined by
$$
c^*{\rightharpoonup}c = c_{(1)}{<}c^*, c_{(2)}{>} \quad \mbox{and}
\quad c{\leftharpoonup}c^* = {<}c^*, c_{(1)}{>}c_{(2)}
$$
for all $c^* \in C^*$ and $c \in C$.

Suppose that $(M, \delta)$ is a left $C$-comodule. For $m \in M$
we use the notation $\delta (m) = m_{(-1)}{\otimes}m_{(0)}$ to
denote $\delta (m) \in C{\otimes}M$. If $(M, \delta)$ is a right
$C$-comodule we denote $\delta (m) \in M{\otimes}C$ by $\delta (m)
= m_{(0)}{\otimes}m_{(1)}$.  Observe that our coproduct and
comodule notations do not conflict.

Bilinear forms play an important role in this paper. We will think
of them in terms of linear forms $\beta : U{\otimes}V
\longrightarrow k$ and will often write $\beta (u, v)$ for $\beta
(u{\otimes}v)$. Note that $\beta$ determines linear maps
$\beta_\ell : U \longrightarrow V^*$ and $\beta_r : V
\longrightarrow U^*$ where $\beta_\ell (u)(v) = \beta (u, v) =
\beta_r(v)(u)$ for all $u \in U$ and $v \in V$. The form $\beta$
is left (respectively right) non-singular if $\beta_\ell$
(respectively $\beta_r$) is one-one and $\beta$ is non-singular if
it is both left and right non-singular.

Suppose that $A$ is an algebra over $k$, $U$ is a right
$A$-module, $V$ is a left $A$-module, and  $\beta : U{\otimes}V
\longrightarrow k$ is a linear form. Then $\beta$ is $A$-balanced
if $\beta(u{\cdot}a, v) = \beta(u, a{\cdot}v)$ for all $u \in U$,
$a \in A$, and $v \in V$.

For subspaces $X \subseteq U$ and $Y\subseteq V$ we define
subspaces $X^\perp \subseteq V$ and $Y^\perp \subseteq U$ by
$$
X^\perp = \{ v \in V \,|\, \beta (X, v) = (0)\,\} \quad \mbox{and}
\quad Y^\perp = \{ u \in U \,|\, \beta (u, Y) (0)\,\}.
$$
Note that there is a form $\overline{\beta} :
U/V^\perp{\otimes}V/U^\perp \longrightarrow k$ uniquely determined
by the commutative diagram
\begin{center}
\begin{picture}(100,70)(0,0)
\put(42,10){\vector(3,1){55}} \put(20,58){\vector(3,-1){77}}
\put(4,45){\vector(0,-1){25}}
\put(-32,5){$U/V^\perp{\otimes}V/U^\perp$}
\put(-10,55){$U{\otimes}V$} \put(102,28){$k$}
\put(50,20){$\overline{\beta}$} \put(50,55){$\beta$}
\end{picture}
\end{center}
where the vertical map is the tensor product of the projections.
Observe that $V^\perp = {\rm Ker}\,\beta_\ell$, $U^\perp = {\rm
Ker}\,\beta_r$, and that $\overline{\beta}$ is non-singular.

Let $A$ be a bialgebra over $k$. A $2$-cocycle for $A$ is a
convolution invertible linear form $\sigma : A{\otimes}A
\longrightarrow k$ which satisfies
$$
\sigma (x_{(1)}, y_{(1)})\sigma(x_{(2)}y_{(2)}, z) =
\sigma(y_{(1)}, z_{(1)})\sigma(x, y_{(2)}z_{(2)})
$$
for all $x, y, z \in A$. If $\sigma$ is a $2$-cocycle for $A$ then
$A^\sigma$ is a bialgebra, where $A^\sigma = A$ as a coalgebra and
multiplication $m^\sigma : A{\otimes}A \longrightarrow A$ is given
by
$$
m^\sigma (x{\otimes}y) = \sigma(x_{(1)},
y_{(1)})x_{(2)}y_{(2)}\sigma^{-1}(x_{(3)}, y_{(3)})
$$
for all $x, y \in A$. See for example \cite{DoiTak}.

Let $U$ and $A$ be bialgebras over $k$ and suppose that $\tau :
U{\otimes}A \longrightarrow k$ is a linear form. Consider the
axioms:
\begin{enumerate}
\item[{\rm (A.1)}] $\tau (u, aa') = \tau (u_{(2)}, a)\tau
(u_{(1)}, a')$ for all $u \in U$ and $a, a' \in A$; \item[{\rm
(A.2)}] $\tau (1, a) = \epsilon(a)$ for all $a \in A$; \item[{\rm
(A.3)}] $\tau (uu', a) = \tau (u, a_{(1)})\tau (u', a_{(2)})$ for
all $u, u'  \in U$ and $a \in A$; \item[{\rm (A.4)}] $\tau (u, 1)
= \epsilon(u)$ for all $u \in U$.
\end{enumerate}
We leave the reader with the exercise of establishing:
\begin{Lemma}\label{LemmaTau}
Let $U$ and $A$ be bialgebras over the field $k$ and suppose $\tau
: U{\otimes}A \longrightarrow k$ is a linear form. Then the
following are equivalent:
\begin{enumerate}
\item[{\rm a)}] (A.1)--(A.4) hold. \item[{\rm b)}] $\tau_\ell (U)
\subseteq A^o$ and $\tau_\ell : U \longrightarrow A^{o\, cop} =
A^{op \, o}$ is a bialgebra map. \item[{\rm c)}] $\tau_r (A)
\subseteq U^o$ and $\tau_r : A \longrightarrow U^{o\, op}$ is a
bialgebra map.
\end{enumerate}
\qed
\end{Lemma}
\medskip

Suppose that (A.1)--(A.4) hold, $\tau$ is convolution invertible,
and define a linear form $\sigma :
(U{\otimes}A){\otimes}(U{\otimes}A) \longrightarrow k$ by $\sigma
(u{\otimes}a, u'{\otimes}a') = \epsilon (a)\tau(u',
a)\epsilon(a')$ for all $u, u' \in U$ and $a, a' \in A$. Then
$\sigma$ is a $2$-cocycle. We denote the $2$-cocycle twist
bialgebra structure on the tensor product bialgebra $U{\otimes}A$
by $H = (U{\otimes}A)^\sigma$. Observe that
\begin{equation}\label{Eq2CocyleMult}
(u{\otimes}a)(u'{\otimes}a') = u\tau(u'_{(1)},
a_{(1)})u'_{(2)}{\otimes}a_{(2)}\tau^{-1}(u'_{(3)}, a_{(3)})a'
\end{equation}
for all $u, u' \in U$ and $a, a' \in A$.

An easy, but important, exercise to do is the following.
\begin{Lemma}\label{LemmaTauInver}
Suppose that $U, A$ are bialgebras over the field $k$ and $\tau :
U{\otimes}A \longrightarrow k$ satisfies (A.1)--(A.4). Then $\tau$
has a convolution inverse if
\begin{enumerate}
\item[{\rm a)}] $U$ is a Hopf algebra with antipode $S$, in which
case $\tau^{-1}(u, a) = \tau(S(u), a)$ for all $u \in U$ and $a
\in A$, or \item[{\rm b)}] $A^{op}$ is a Hopf algebra with
antipode $T$, in which case $\tau^{-1}(u, a) = \tau(u, T(a))$ for
all $u \in U$ and $a \in A$.
\end{enumerate}
\qed
\end{Lemma}
\medskip

The quantum double provides an important example of a $2$-cocycle
twist bialgebra \cite{DoiTak}.
\begin{Ex}\label{ExDoubleTwist}
Let $A$ be a finite-dimensional Hopf algebra over $k$, let $U =
A^{o\, cop}$, and let $\tau : U{\otimes}A \longrightarrow k$ be
defined by $\tau (p, a) = p(a)$ for all $p \in U$ and $a \in A$.
Then $\tau_\ell : U \longrightarrow A^{o \, cop}$ is the identity
map and $(U{\otimes}A)^\sigma = D(A)$.
\end{Ex}

Observe that finite-dimensionality was not necessary to define a
bialgebra structure on $D(A) = A^{o \, cop}{\otimes}A$. For any
Hopf algebra $A$ with bijective antipode we let $D(A) =
(U{\otimes}A)^\sigma$, where $\tau$ is defined as above.

Suppose that $U$, $\overline{U}$ and $A$, $\overline{A}$ are
algebras over $k$. Suppose further that $\tau :U{\otimes}A
\longrightarrow k$ and $\overline{\tau} :
\overline{U}{\otimes}\overline{A} \longrightarrow k$ are
convolution invertible linear forms satisfying (A.1)--(A.4). Set
$H = (U{\otimes}A)^{\sigma}$ and $\overline{H} =
(\overline{U}{\otimes}\overline{A})^{\overline{\sigma}}$. Suppose
that $f : U \longrightarrow \overline{U}$ and $g : A
\longrightarrow \overline{A}$ are bialgebra maps such that
$\overline{\tau}(f(u), g(a)) = \tau (u, a)$ for all $u \in U$ and
$a \in A$. Then $f{\otimes}g : H \longrightarrow \overline{H}$ is
a bialgebra map.

As a consequence $f: H \longrightarrow D(A)$ defined by
$f(u{\otimes}a) = \tau_\ell (u){\otimes}a$ for all $u \in U$ and
$a \in A$ is a bialgebra map. In this paper we are interested in
left modules over $H$. A good source is modules for the double in
light of the map $f$.
\section{Algebra Structures on the Vector Space $U{\otimes}A$, where $U$ and $A$ are
algebras over $k$}\label{SecUOtimesA} Suppose that $U$ and $A$ are
algebras over the field $k$. In this section we are interested in
algebra structures $H = U{\otimes}A$ on the tensor product of
their underlying vector spaces which satisfy
\begin{equation}\label{EqUandAInclusions}
(u{\otimes}a)(u'{\otimes}a') = uu'{\otimes}aa' \qquad
\mbox{whenever} \qquad \mbox{$a = 1$ or $u' = 1$.}
\end{equation}
For such an algebra the maps $U \longrightarrow H$ and $A
\longrightarrow H$ given by $u \mapsto u{\otimes}1$ and $a \mapsto
1{\otimes}a$ respectively are algebra maps. As a consequence $H$
is a left $U$-module and a right $A$-module by pullback action;
thus
$$
u{\cdot}(u'{\otimes}a) = (u{\otimes}1)(u'{\otimes}a) =
uu'{\otimes}a
$$
and
$$
(u{\otimes}a){\cdot}a' = (u{\otimes}a)(1{\otimes}a') =
u{\otimes}aa'
$$
for all $u, u' \in U$ and $a, a' \in A$.

We list several examples of these algebras.
\begin{Ex}\label{UTensorA}
Let $U$ and $A$ be algebras over the field $k$. Then the tensor
product algebra structure on $H = U{\otimes}A$ satisfies
(\ref{EqUandAInclusions}).
\end{Ex}
\begin{Ex}\label{HOp}
Let $U$ and $A$ be algebras over the field $k$ and $H =
U{\otimes}A$ be an algebra structure on the tensor product of
their underlying vector spaces which satisfies
(\ref{EqUandAInclusions}). Endow the vector space
$A^{op}{\otimes}U^{op} = A{\otimes}U$ with the unique algebra
structure which makes the twist map $\tau_{A, U} : A{\otimes}U
\longrightarrow (U{\otimes}A)^{op}$ an algebra isomorphism. This
algebra $H^{\widetilde{op}} = A^{op}{\otimes}U^{op}$ satisfies
(\ref{EqUandAInclusions}).
\end{Ex}
\begin{Ex}\label{ExDouble}
Let $A$ be a finite-dimensional Hopf algebra with antipode $s$
over $k$ and let $U = A^{*\, cop}$. As a vector space the Drinfeld
double is $D(A) = U{\otimes}A$ and its product is determined by
$$
(p{\otimes}a)(q{\otimes}b) =
p(a_{(1)}{\cdot}q{\cdot}s^{-1}(a_{(3)})){\otimes}a_{(2)}b
$$
for all $p, q \in U$ and $a, b \in A$. Thus the underlying algebra
structure of $D(A)$ satisfies (\ref{EqUandAInclusions}).
\end{Ex}

The next example is the most important one for us. By virtue of
Example \ref{ExDoubleTwist} the preceding example is a special
case of it.
\begin{Ex}\label{ExTwoCocycle}
Let $U$ and $A$ be bialgebras over the field $k$ and suppose that
$\tau : U{\otimes}A \longrightarrow k$ is a convolution invertible
map such (A.1)--(A.4) are satisfied. Then $H =
(U{\otimes}A)^\sigma$ defined in Section \ref{SecPrelim} is a
bialgebra whose underlying algebra structure satisfies
(\ref{EqUandAInclusions}).
\end{Ex}
\begin{Ex}\label{ExBiProd}
Let $H$ be a Hopf algebra with bijective antipode over $k$ and
suppose that $R \in {}_H^H{\mathcal YD}$ is a bialgebra in the
Yetter--Drinfeld category. Then the bi-product $A = R{\#}H$
satisfies (\ref{EqUandAInclusions}). More generally smash products
satisfy (\ref{EqUandAInclusions}).
\end{Ex}
\noindent A basic reference for the bi-product is \cite{RP}. See
\cite{ASSurvey} for a discussion of the Yetter--Drinfeld category
${}_H^H{\mathcal YD}$ and bi-products.

Suppose that $U$, $\overline{U}$, $A$, and $\overline{A}$ are
algebras over $k$ and that the vector spaces $U{\otimes}A$ and
$\overline{U}{\otimes}\overline{A}$ have algebra structures which
satisfy (\ref{EqUandAInclusions}). A morphism $F : U{\otimes}A
\longrightarrow \overline{U}{\otimes}\overline{A}$ of these
algebras is a map of algebras which satisfies $F(U{\otimes}1)
\subseteq \overline{U}{\otimes}1$ and $F(1{\otimes}A) \subseteq
1{\otimes}\overline{A}$. By virtue of (\ref{EqUandAInclusions}) a
morphism $F : U{\otimes}A \longrightarrow
\overline{U}{\otimes}\overline{A}$ has the form $F = f{\otimes}g$,
where $f : U \longrightarrow \overline{U}$ and $g : A
\longrightarrow \overline{A}$ are algebra maps. Conversely, if $f
: U \longrightarrow \overline{U}$ and $g : A \longrightarrow
\overline{A}$ are algebra maps, and the function $f{\otimes}g :
U{\otimes}A \longrightarrow \overline{U}{\otimes}\overline{A}$ is
an algebra map, then $F = f{\otimes}g$ is a morphism.

Suppose that $U'$ is a subalgebra of $U$ and $A'$ is a subalgebra
of $A$ such that $U'{\otimes}A'$ is a subalgebra of $U{\otimes}A$.
Then $U'{\otimes}A'$ satisfies (\ref{EqUandAInclusions}) and the
tensor product of the inclusion maps $i_{U'}{\otimes}i_{A'} :
U'{\otimes}A' \longrightarrow U{\otimes}A$ is a morphism.

Now suppose $I$ is an ideal of $U$, that $J$ is an ideal of $A$,
and $K = I{\otimes}A + U{\otimes}J$ is an ideal of $U{\otimes}A$.
Let $\pi_I : U \longrightarrow U/I$ and $\pi_J : A \longrightarrow
A/J$ be the projections. Endow $(U/I){\otimes} (A/J)$ with the
algebra structure which makes the linear isomorphism
$(U{\otimes}A)/K \longrightarrow (U/I){\otimes}(A/J)$ given by
$u{\otimes}a + K \mapsto (u + I){\otimes}(a + J)$ an isomorphism
of algebras. Then $(U/I){\otimes}(A/J)$ satisfies
(\ref{EqUandAInclusions}) and the tensor product of projections
$\pi_I{\otimes}\pi_J : U{\otimes}A \longrightarrow
(U/I){\otimes}(A/J)$ is a morphism.

The constructions of the preceding paragraph can be combined to
give a first isomorphism theorem for the morphisms of this
section. Let $F : U{\otimes}A \longrightarrow
\overline{U}{\otimes}\overline{A}$ be a morphism and write $F =
f{\otimes}g$, where $f : U \longrightarrow \overline{U}$ and $g :
A \longrightarrow \overline{A}$ are algebra maps. Then ${\rm
Im}\,f$ and ${\rm Im}\,g$ are subalgebras of $U$ and $A$
respectively and ${\rm Im}\,F = {\rm Im}\,f{\otimes}{\rm Im}\,g$
is a subalgebra of $\overline{U}{\otimes}\overline{A}$ which thus
satisfies (\ref{EqUandAInclusions}). Now ${\rm Ker}\,$ and ${\rm
Ker}\,g$ are ideals of $U$ and $A$ respectively and ${\rm Ker}\,F
= {\rm Ker}\,f{\otimes}A + U{\otimes}{\rm Ker}\,g$ is an ideal of
$U{\otimes}A$. Identifying $(U{\otimes}A)/K$ and $(U/{\rm
Ker}\,f){\otimes}(A/{\rm Ker}\,g)$ as above, note there is a
unique morphism
$$
\overline{F} : (U/{\rm Ker}\,f){\otimes}(A/{\rm
Ker}\,g) \longrightarrow {\rm Im}\,F = {\rm Im}\,f{\otimes}{\rm
Im}\,g
$$
such that $\overline{F}{\circ}(\pi_{{\rm Ker}\,
f}{\otimes}\pi_{{\rm Ker}\, g}) = F$.
\subsection{The Construction of ${\rm L}(\rho, \chi)$ and $R(\chi, \rho)$}\label{SubSecPairs}
Let $U$ and $A$ be algebras over $k$ and $H = U{\otimes}A$ be an
algebra structure on the tensor product of their underlying vector
spaces which satisfies (\ref{EqUandAInclusions}). We will identify
$U$ and $A$ with their images under the algebra maps $U
\longrightarrow U{\otimes}A$ and $A \longrightarrow U{\otimes}A$
given by $u \mapsto u{\otimes}1$ and $a \mapsto 1{\otimes}a$
respectively. In this section we construct special representations
of $H$ by induction on one-dimensional representations of $A$ and
of $U$ determined by pairs $(\rho, \chi)$, where $\rho$ is a
character of $U$ and $\chi$ is a character of $A$.

Suppose that $\rho \in \Alg(U, k)$ and $\chi \in \Alg(A, k)$. We
give $k$ the left $A$-module structure $(k,{\cdot}_\chi)$ where
$$
a{\cdot}_\chi1 = \chi (a)1 \qquad \mbox{for all $a \in A$.}
$$
Recall that $H = U{\otimes}A$ is a right $A$-module via
$(u{\otimes}a){\cdot}a' = u{\otimes}aa'$ for all $u \in U$ and $a,
a' \in A$. We identify the left $H$-module $H{\otimes}_Ak$ and $U$
by the linear isomorphism $H{\otimes}_Ak \longrightarrow U$ given
by $(u{\otimes}a)_A{\otimes}1 \mapsto u\chi (a)$ for all $u \in U$
and $a \in A$ and denote the resulting right $H$-module structure
on $U$ by $(U,  {\cdot}_\chi)$. We write $U_\chi$ for $U$ with
this module action implicitly understood. Note that
\begin{equation}\label{EqUChiAction}
(u{\otimes}a){\cdot}_\chi u' = u\!\left(\PH ({\rm
I}_U{\otimes}\chi )((1{\otimes}a)(u'{\otimes}1))\right),
\end{equation}
and consequently
\begin{equation}\label{EqUandAChiTrans}
u{\cdot}_\chi u' = uu' \qquad \mbox{and} \qquad a{\cdot}_\chi 1 =
\chi (a)1,
\end{equation}
for all $u, u' \in U$ and $a \in A$. Thus $1$ generates $U_\chi$
as a left $H$-module and $k1$ is a one-dimensional left
$A$-submodule of $U_\chi$.

Let $I(\rho, \chi)$ be the sum of all the left $H$-submodules of
$U_\chi$ contained in ${\rm Ker}\,\rho$, let
$$
L(\rho, \chi) = U_\chi/I(\rho, \chi)
$$
be the resulting quotient left $H$-module, and let $\pi_{(\rho,
\chi)} : U_\chi \longrightarrow L(\rho, \chi)$ be the projection.
Using (\ref{EqUandAChiTrans}) we see that $M = L(\rho, \chi)$ is a
cyclic left $H$-module generated by $\pi_{(\rho, \chi)} (1)$, that
$k\pi_{(\rho, \chi)}  (1)$ is a one-dimensional left $A$-submodule
of $M$, and that $N = \pi_{(\rho, \chi)} ({\rm Ker}\,\rho)$ is a
codimension one left $U$-submodule of $M$ with the property that
the only left $H$-submodule of $M$ contained in $N$ is $(0)$.

In a similar manner we define a right $H$-module structure on $A$.
Regard $k$ as the right $U$-module $(k, {\cdot}_\rho)$ where
$$
1_\rho {\cdot}u = \rho (u)1 \qquad \mbox{for all $u \in U$.}
$$
Recall that $H = U{\otimes}A$ is a left $U$-module via
$u{\cdot}(u'{\otimes}a)= uu'{\otimes}a$ for all $u, u' \in U$ and
$a \in A$. We identify the right $H$-module $k{\otimes}_UH$ and
the vector space $A$ by the linear isomorphism $k{\otimes}_UH
\longrightarrow A$ given by $1{\otimes}_U(u{\otimes}a) \mapsto
\rho(u)a$ for all $u \in U$ and $a \in A$ and we denote the
resulting module structure on $A$ by $(A, {\cdot}_\rho)$. Observe
that
\begin{equation}\label{EqARhoAction}
a{\cdot}_\rho(u{\otimes}a') = \left(\PH (\rho{\otimes}{\rm
I}_A)((1{\otimes}a)(u{\otimes}1))\right)\!a',
\end{equation}
and thus
\begin{equation}\label{EqAandURhoTrans}
a{\cdot}_\rho a' = aa' \qquad \mbox{and} \qquad 1{\cdot}_\rho u =
\rho (u)1,
\end{equation}
for all $a, a' \in A$ and $u \in U$. As a consequence $1$
generates $A_\rho$ as a right $H$-module and $k1$ is a
one-dimensional left $U$-submodule of $A_\rho$.

Let $J(\chi, \rho)$ be the sum of all the right $H$-submodules of
$A_\rho$ contained in ${\rm Ker}\,\chi$, let
$$
R(\chi, \rho) = A_\rho/J(\chi, \rho)
$$
be the quotient right $H$-module, and $\pi_{(\chi, \rho)} : A_\rho
\longrightarrow R(\chi, \rho)$ be the projection. Using
(\ref{EqAandURhoTrans}) we see that $M = R(\chi, \rho)$ is a
cyclic right $H$-module generated by $\pi_{(\chi, \rho)} (1)$,
that $k\pi_{(\chi, \rho)}  (1)$ is a one-dimensional right
$U$-submodule of $M$, and that $N = \pi_{(\chi, \rho)} ({\rm
Ker}\,\chi)$ is a codimension one left $A$-submodule of $M$ with
the property that the only right $H$-submodule of $M$ contained in
$N$ is $(0)$.
\subsection{A Bilinear Form Arising from Character Pairs $(\rho, \chi)$}\label{SubSecBilinear}
We continue with the notation and assumptions of the preceding
section. Let $\chi \in \Alg(A, k)$, let $\rho \in \Alg \,(U, k)$,
and let $\Psi : A{\otimes} U \longrightarrow k$ be the linear form
defined by
$$
\Psi (a, u) = (\rho{\otimes}\chi)\left(\PH
(1{\otimes}a)(u{\otimes}1)\right)
$$
for all $a \in A$ and $u \in U$. The $L(\rho, \chi)$ and $R(\chi,
\rho)$ constructions of the preceding section are related in
fundamental ways through this form.

Regard $U^*$ as a right $H$-module with the transpose action
determined by $U_\chi$ and likewise regard $A^*$ as left
$H$-module with the transpose action determined by $A_\rho$. Our
notations for these actions are given in the equations
$$
(u^*{\cdot}_\chi h)(u) = u^*(h{\cdot}_\chi u) \qquad \mbox{and}
\qquad (h{\cdot}_\rho a^*)(a) = a^* (a{\cdot}_\rho h)
$$
for all $u^* \in U^*$, $h \in H$, $u \in U$, $a^* \in A^*$, and $a
\in A$.

\begin{Prop}\label{PropBilinear}
Let $U$ and $A$ be algebras over the field $k$ and let $H =
U{\otimes}A$ be an algebra structure on the tensor product of
their underlying vector spaces which satisfies
(\ref{EqUandAInclusions}). Let $\chi \in \Alg(A, k)$, let $\rho
\in \Alg(U, k)$, and let $\Psi : A{\otimes}U \longrightarrow k$ be
the bilinear form defined above. Then:
\begin{enumerate}
\item[{\rm a)}] $\Psi (a, u) = \rho((1{\otimes}a){\cdot}_\chi u) =
\chi (a{\cdot}_\rho (u{\otimes}1))$ for all $a \in A$ and $u \in
U$. \item[{\rm b)}] $\Psi(a{\cdot}_\rho h, u) = \Psi(a,
h{\cdot}_\chi u)$ for all $a \in A$, $h \in H$, and $u \in U$;
that is $\Psi$ is $H$-balanced. \item[{\rm c)}] $A^\perp I(\rho,
\chi)$ and thus $U/A^\perp = L(\rho, \chi)$. Furthermore there is
an isomorphism of left $H$-modules $L(\rho, \chi) \longrightarrow
U{\cdot}_\rho\chi$ given by $u + I(\rho, \chi) \mapsto \Psi_r(u) =
u{\cdot}_\rho\chi$ for all $u \in U$. \item[{\rm d)}] $U^\perp =
J(\chi, \rho)$ and thus $A/U^\perp = R(\chi, \rho)$. Furthermore
there is an isomorphism of right $H$-modules $R(\chi, \rho)
\longrightarrow \rho{\cdot}_\chi A$ given by $a + J(\chi, \rho)
\mapsto \Psi_\ell (a) = \rho{\cdot}_\chi a$ for all $a \in A$.
\item[{\rm e)}] There is a non-singular $H$-balanced bilinear form
$$
\overline{\Psi} : R(\chi, \rho){\otimes} L(\rho,
\chi)\longrightarrow k
$$ determined by the commutative diagram
\begin{center}
\begin{picture}(100,70)(0,0)
\put(52,12){\vector(3,1){45}} \put(20,58){\vector(3,-1){77}}
\put(4,45){\vector(0,-1){25}}
\put(-34, 5){$R(\chi, \rho){\otimes} L(\rho, \chi)$}
\put(-10,55){$A{\otimes}U$} \put(102,26){$k$}
\put(60,20){$\overline{\Psi}$} \put(55,50){$\Psi$}
\end{picture}
\end{center}
where the vertical arrow is the tensor product of projections.
\end{enumerate}
\end{Prop}

\pf Part a) follows by the definition of the form and of the
module actions. To show part b) we use several standard
isomorphisms involving tensor products. Note that $k$ is a left
and right $A$-module via $a{\cdot}_\chi 1 = 1{\cdot}_\chi a = \chi
(a)1$ for all $a \in A$ and $k$ is also a left and right
$U$-module via $u{\cdot}_\rho 1 = 1{\cdot}_\rho u = \rho (u)1$ for
all $u \in U$. Consider the composites $f : A{\otimes}U
\longrightarrow k{\otimes}_U(H{\otimes}_Ak)$ and $g :
k{\otimes}_U(H{\otimes}_Ak) \longrightarrow k$ defined by
\begin{eqnarray*}
\lefteqn{A{\otimes}U \simeq
(k{\otimes}_UH){\otimes}(H{\otimes}_Ak)
\longrightarrow (k{\otimes}_UH) {\otimes}_H(H{\otimes}_Ak)}  \\
& & \qquad \simeq k{\otimes}_U(H{\otimes}_H(H{\otimes}_Ak)) \simeq
k{\otimes}_U(H{\otimes}_Ak)
\end{eqnarray*}
and
$$
k{\otimes}_U(H{\otimes}_Ak) \stackrel{{\rm I}_k
{\otimes}((\chi{\otimes}\rho){\otimes}{\rm I}_k)}{\longrightarrow}
k{\otimes}_U(k{\otimes}_Ak) \stackrel{{\rm
I}_k{\otimes}m}{\longrightarrow} k{\otimes}_Uk
\stackrel{m}{\longrightarrow} k
$$
respectively, where $m$ is multiplication. Since
$$
(g{\circ}f)(a{\otimes}u) =
g(1{\otimes}_U((1{\otimes}a)(u{\otimes}1)_A{\otimes}1)) =
(\rho{\otimes}\chi)((1{\otimes}a)(u{\otimes}1))
$$
for all $a \in A$ and $u \in U$ we conclude that $\Psi =
g{\circ}f$. From the second isomorphism in the definition of $f$
we see that $f((a{\cdot}_\rho h){\otimes}u) =
f(a{\otimes}(h{\cdot}_\chi u))$ for all $a \in A$ and $u \in U$.
Part b) now follows.

We next show part c). Let $u \in U$. Since $\Psi_r(u)(1) =
\rho(u)$ by part a) it follows that ${\rm Ker}\,\Psi_r \subseteq
{\rm Ker}\,\rho$. Now ${\rm Ker}\,\Psi_r$ is a left $H$-submodule
of $U_\chi$ by part b). Suppose that $L$ is a left $H$-submodule
of $U_\chi$ and $L \subseteq {\rm Ker}\,\rho$. Then by part a)
again $\Psi(A, L) \subseteq \rho ((1{\otimes}A){\cdot}_\chi L)
\subseteq \rho (L) = (0)$. Therefore $L \subseteq {\rm
Ker}\,\Psi_r$. We have shown that $A^\perp = {\rm Ker}\,\Psi_r =
I(\rho, \chi)$.

To complete the proof of part c) we note that $\Psi_r : U_\chi
\longrightarrow A^*$ is a map of left $H$-modules by part b). Let
$u \in U$. Since $\Psi_r(u) = \Psi_r(u{\cdot}_\chi 1) =
u{\cdot}_\rho \Psi_r(1) = u{\cdot}_\rho \chi$, the isomorphism of
left $H$-modules
$$
L(\rho, \chi) = U_\chi/A^\perp = U/{\rm Ker}\,\Psi_r \simeq
U{\cdot}_\rho \chi
$$
is given by $u + I(\rho, \chi) \mapsto u{\cdot}_\rho \chi$. We
have established part c). The proof of part d) is similar. Part e)
follows from parts b)--d). \qed
\medskip

By virtue of part c) the quotient $L(\rho, \chi)$ can be realized
as submodule.
\subsection{The Connection Between Morphisms of the Algebra Structures and the
Constructions  ${\rm L}(\rho, \chi)$, $R(\chi,
\rho)$}\label{SubSecMorLR} Let $U$, $\overline{U}$, $A$, and
$\overline{A}$ be algebras over $k$  such that $H = U{\otimes}A$
and $\overline{H} = \overline{U}{\otimes}\overline{A}$ are
algebras which satisfy (\ref{EqUandAInclusions}).  Suppose that $F
: U{\otimes}A \longrightarrow \overline{U}{\otimes}\overline{A}$
is a morphism and let $f : U \longrightarrow \overline{U}$ and $g
: A \longrightarrow \overline{A}$ be the unique algebra maps which
satisfy $F = f{\otimes}g$.

Now suppose that $\overline{\rho} \in \Alg (\overline{U}, k)$ and
$\overline{\chi} \in \Alg (\overline{A}, k)$. Then $\rho =
\overline{\rho}{\circ}f \in \Alg (U, k)$ and $\chi =
\overline{\chi}{\circ}g \in \Alg (A, k)$. In this section we
examine the connection between  $L(\rho, \chi)$ and
$L(\overline{\rho}, \overline{\chi})$ and also the relationship
between $R(\chi, \rho)$ and $R(\overline{\rho}, \overline{\chi})$.

First we consider $U_\chi$ and $\overline{U}_{\overline{\chi}}$.
Let $u, u' \in U$ and $a \in A$. By the calculation
\begin{eqnarray*}
f((u{\otimes}a){\cdot}_\chi u') & = &
f(u\!\left(\PH (I_U{\otimes}\chi )((1{\otimes}a)(u'{\otimes}1))\right)) \\
& = & f(u\!\left(\PH (I_U{\otimes}(\overline{\chi}{\circ}g)
)((1{\otimes}a)
(u'{\otimes}1))\right)) \\
& = & f(u)\!\left((\PH I_U{\otimes}\overline{\chi})\left(\PH
(f{\otimes}g)
((1{\otimes}a)(u'{\otimes}1))\right)\right) \\
& = & f(u)\!\left((\PH I_U{\otimes}\overline{\chi})\left(\PH F
((1{\otimes}a)
(u'{\otimes}1))\right)\right) \\
& = & f(u)\!\left((\PH I_U{\otimes}\overline{\chi})\left(\PH F
(1{\otimes}a)
F(u'{\otimes}1)\right)\right) \\
& = & f(u)\!\left((\PH I_U{\otimes}\overline{\chi})\left(\PH
(1{\otimes}f(a))
(g(u'){\otimes}1)\right)\right) \\
& = & F(u{\otimes}a){\cdot}_{\overline{\chi}}f(u')
\end{eqnarray*}
we see that
\begin{equation}\label{EqPullBackf}
f((u{\otimes}a){\cdot}_\chi u') =
F(u{\otimes}a){\cdot}_{\overline{\chi}}f(u').
\end{equation}
Thus $f$ is $F$-linear and the preceding equation has a simple
interpretation in terms of module maps. Regard
$\overline{U}_{\overline{\chi}}$ as a left $H$-module by pullback
along $F$. Then $f : U_\chi \longrightarrow
\overline{U}_{\overline{\chi}}$ is a map of left $H$-modules.

Since $\overline{\rho}{\circ}f = \rho$, we have $f({\rm
Ker}\,\rho) \subseteq {\rm Ker}\,\overline{\rho}$. Thus $f(I(\rho,
\chi)) \subseteq I(\overline{\rho}, \overline{\chi})$ as $f$ is
$F$-linear. Consequently $f$ gives rise to an $F$-linear map $L(f)
: L(\rho, \chi) \longrightarrow L(\overline{\rho},
\overline{\chi})$.
\begin{Prop}\label{PropfLift}
Suppose $U$, $\overline{U}$, $A$, and $\overline{A}$ are algebras
over $k$. Let $H = U{\otimes}A$ and $\overline{H} =
\overline{U}{\otimes}\overline{A}$ be algebra structures on the
tensor product of underlying vector spaces which satisfy
(\ref{EqUandAInclusions}).  Suppose further that $F : U{\otimes}A
\longrightarrow \overline{U}{\otimes}\overline{A}$ is a morphism
and let $f : U \longrightarrow \overline{U}$ and $g : A
\longrightarrow \overline{A}$ be the unique algebra maps which
satisfy $F = f{\otimes}g$. Then:
\begin{enumerate}
\item[{\rm a)}] $f : U_\chi \longrightarrow
\overline{U}_{\overline{\chi}}$ is $F$-linear. \item[{\rm b)}]
There is a unique $F$-linear map $L(f) : L(\rho, \chi)
\longrightarrow L(\overline{\rho}, \overline{\chi})$ which makes
the diagram
\begin{center}
\begin{picture}(100,70)(0,0)
\put(15,0){$U_\chi$} \put(0,50){$L(\rho, \chi)$}
\put(100,0){$\overline{U}_{\overline{\chi}}$}
\put(85,50){$L(\overline{\rho}, \overline{\chi})$}
\put(40,3){\vector(1,0){50}} \put(45,53){\vector(1,0){35}}
\put(20,15){\vector(0,1){25}} \put(105,15){\vector(0,1){25}}
\put(60,8){$f$} \put(50,58){$L(f)$}
\end{picture}
\end{center}
commute, where the vertical arrows are the projection maps.
\item[{\rm c)}] Suppose that $f$ is onto. Then $L(f)$ is an
isomorphism.
\end{enumerate}
\end{Prop}

\pf We have established parts a) and b) in the discussion
preceding the statement of the proposition. It remains to show
part c).

We first observe that $\overline{\rho}{\circ}f = \rho$ implies
${\rm Ker}\,\rho = f^{-1}({\rm Ker}\,\overline{\rho})$. Therefore
$f^{-1}(I(\overline{\rho}, \overline{\chi}))$ is a left
$H$-submodule of $U_\chi$ contained in ${\rm Ker}\,\rho$ which
implies $f^{-1}(I(\overline{\rho}, \overline{\chi})) \subseteq
I(\rho, \chi)$. We have seen that $f(I(\rho, \chi)) \subseteq
I(\overline{\rho}, \overline{\chi})$ in any event. Therefore
$f^{-1}(I(\overline{\rho}, \overline{\chi})) = I(\rho, \chi)$
which implies that $L(f)$ is an isomorphism. \qed
\medskip

Suppose that $F$ is onto. Then the hypothesis of part c) is met.
Regard $I(\overline{\rho}, \overline{\chi})$ as a left $H$-module
by pullback along $F$. Then the $H$-submodules and
$\overline{H}$-submodules of $L(\overline{\rho}, \overline{\chi})$
are the same and the $H$-module $L(\rho, \chi)$ is understood in
terms of the left $\overline{H}$-module $L(\overline{\rho},
\overline{\chi})$ and the algebra map $F$.

There is an analog of the preceding proposition for $R(\chi,
\rho)$ and $R(\overline{\chi}, \overline{\rho})$. One can show
that
\begin{equation}\label{EqPullBackg}
g(a{\cdot}_\rho (u{\otimes}a')) =
g(a){\cdot}_{\overline{\rho}}F(u{\otimes}a')
\end{equation}
for all $a, a' \in A$ and $u \in U$ by mimicking the calculation
which establishes (\ref{EqPullBackf}).  By modifying the proof of
the preceding proposition one can easily show:
\begin{Prop}\label{PropgLift}
Suppose $U$, $\overline{U}$, $A$, and $\overline{A}$ are algebras
over $k$. Let $H = U{\otimes}A$ and $\overline{H} =
\overline{U}{\otimes}\overline{A}$ be algebra structures on the
tensor product of underlying vector spaces which satisfy
(\ref{EqUandAInclusions}).  Suppose further that $F : U{\otimes}A
\longrightarrow \overline{U}{\otimes}\overline{A}$ is a morphism
and let $f : U \longrightarrow \overline{U}$ and $g : A
\longrightarrow \overline{A}$ be algebra maps which satisfy $F =
f{\otimes}g$. Then:
\begin{enumerate}
\item[{\rm a)}] $g : A_\rho \longrightarrow
\overline{A}_{\overline{\rho}}$ is $F$-linear. \item[{\rm b)}]
There is a unique $F$-linear map $R(g) : R(\chi, \rho)
\longrightarrow R(\overline{\chi}, \overline{\rho})$ which makes
the diagram
\begin{center}
\begin{picture}(100,70)(0,0)
\put(15,0){$A_\rho$} \put(0,50){$R(\chi, \rho)$}
\put(100,0){$\overline{A}_{\overline{\rho}}$}
\put(85,50){$R(\overline{\chi}, \overline{\rho})$}
\put(40,3){\vector(1,0){50}} \put(45,53){\vector(1,0){35}}
\put(20,15){\vector(0,1){25}} \put(105,15){\vector(0,1){25}}
\put(60,8){$g$} \put(50,58){$R(g)$}
\end{picture}
\end{center}
commute, where the vertical arrows are the projection maps.
\item[{\rm c)}] Suppose that $g$ is onto. Then $R(g)$ is an
isomorphism.
\end{enumerate}
\qed
\end{Prop}
\medskip

We conclude this section by noting the relationship between the
linear forms $\Psi : A{\otimes}U \longrightarrow k$ and
$\overline{\Psi} : \overline{U}{\otimes}\overline{A}
\longrightarrow k$ defined by
$$
\Psi (a, u) = (\rho{\otimes}\chi)((1{\otimes}a)(u{\otimes}1))
\quad \mbox{and} \quad \overline{\Psi} (\overline{a},
\overline{u}) =
(\overline{\rho}{\otimes}\overline{\chi})((1{\otimes}\overline{a})(\overline{u}{\otimes}1))
$$
for all $a \in A$, $u \in U$, $\overline{a} \in \overline{A}$, and
$\overline{u} \in \overline{U}$. Since $F = f{\otimes}g$ is an
algebra map, the calculation
\begin{eqnarray*}
(f{\otimes}g)((1{\otimes}a)(u{\otimes}1)) & = &
((f{\otimes}g)(1{\otimes}a))((f{\otimes}g)(u{\otimes}1)) \\
& = & (1{\otimes}g(a))(f(u){\otimes}1)
\end{eqnarray*}
shows that
$$
\Psi (a, u) = \overline{\Psi}(g(a), f(u)) = \overline{\Psi}(F(a),
F(u))
$$
for all $a \in A$ and $u \in U$. In the last expression we regard
$a, u \in H$ by the identifications $a = 1{\otimes}a$ and $u =
u{\otimes}1$.
\subsection{The Categories ${}_H\underline{\mathcal M}$ and $\underline{\mathcal M}_H$
and Duality}\label{SubSecCatM} Let $H = U{\otimes}A$ be an algebra
structure defined on the tensor product of their underlying vector
spaces and suppose that (\ref{EqUandAInclusions})  holds for $H$.
Let ${}_H\underline{\mathcal M}$ be the category whose objects $M$
are left $H$-modules which are generated by a one-dimensional left
$A$-submodule $km$ and have a codimension-one left $U$-submodule
$N$ with the property that $(0)$ is the only left $H$-submodule of
$M$ contained in $N$. We take maps of left $H$-modules to be our
morphisms. The category $\underline{\mathcal M}_H$ is defined in
the same manner with ``right" replacing ``left" and with the roles
of $U$ and $A$ reversed.

For $M$ as described above observe that
\begin{equation}\label{EqChiRho}
{\rm ann}_A(km) = {\rm Ker}\,\chi \qquad \mbox{and} \qquad {\rm
ann}_U(M/N) = {\rm Ker}\,\rho
\end{equation}
for some characters $\chi \in \Alg(A, k)$ and $\rho \in \Alg(U,
k)$. For any pair of characters $\chi \in \Alg \,(A, k)$ and $\rho
\in \Alg(U, k)$ observe that $M = L(\rho, \chi)$ is an object of
${}_H{\mathcal M}$ which satisfies (\ref{EqChiRho}).

Conversely, suppose that $M$ is an object of
${}_H\underline{\mathcal M}$ which satisfies (\ref{EqChiRho}).
Then the rule $H{\otimes}_Ak \longrightarrow M$ given by
$h{\otimes}_A1 \mapsto h{\cdot}m$ is a well-defined map of left
$H$-modules and the composite $f : U \simeq H{\otimes}_A k
\longrightarrow M$, which is given by $f(u) = u{\cdot}m$ for all
$u \in U$, has kernel ${\rm I}(\chi, \rho)$. Therefore $f$ lifts
to an isomorphism of left $H$-modules $L(\rho, \chi)  \simeq M  $.
Observe that $f^{-1}(N) = {\rm Ker}\,\rho$. Therefore $N$ is the
only codimension one left $U$-submodule $N'$ of $M$ such that
${\rm ann}_U(M/N') = {\rm Ker}\,\rho$.

To discuss duality we need to specify a particular one-dimensional
left $A$-submodule and a particular codimension-one $U$-submodule
of each object $M$ of ${}_H\underline{\mathcal M}$. Let
${}_H\underline{\mathcal M}'$ be the category whose objects are
triples $(M, km, N)$, where $M$ is a left $H$-module, $km$ is a
left $A$-submodule of $M$ which generates $M$ as a left
$H$-module, and $N$ is a codimension-one left $U$-submodule of $M$
such that $(0)$ is the only left $H$-submodule of $M$ contained in
$N$. A morphism $f : (M, km, N) \longrightarrow (M', km', N')$ of
${}_H\underline{\mathcal M}'$ is a map of left $H$-modules which
satisfies $f(km) \subseteq km'$ and $f(N) \subseteq N'$. We define
a category $\underline{\mathcal M}'_H$ in the same manner
replacing ``left" by ``right".

There is a natural contravariant functor ${}_H\underline{\mathcal
M}' \longrightarrow \underline{\mathcal M}'_H$. To describe it we
start in a slightly more general context.

Consider a triple $(M, km, N)$, where $M$ is a cyclic left
$H$-module generated by $m$, where $N$ is a codimension one left
$U$-submodule of $M$, and (\ref{EqChiRho}) is satisfied for some
$\chi \in \Alg(A, k)$ and $\rho \in \Alg(U, k)$. Thus we are not
requiring that the only left $H$-submodule of $M$ contained in $N$
is $(0)$. We regard $M^* \in {\mathcal M}_H$ by the transpose
action on $M \in {}_H{\mathcal M}$. Since $N$ is a subspace of $M$
of codimension one and $m \not \in N$ there is a non-zero
$m^{\bullet} \in M^*$ uniquely determined by $m^{\bullet}(N) =
(0)$ and $m^{\bullet}(m) = 1$.

Consider the right $H$-submodule $M^{\bullet} =
m^{\bullet}{\cdot}H$ of $M^*$. Now ${\rm Ker}\,\rho = {\rm
ann}_U(M/N)$ implies that $({\rm Ker}\,\rho){\cdot}M \subseteq N$.
From the calculation
$$
(m^{\bullet}{\cdot}({\rm Ker}\,\rho))(M) = m^{\bullet}(({\rm
Ker}\,\rho){\cdot}M) \subseteq m^{\bullet}(N) = (0)
$$
we conclude that $m^{\bullet}{\cdot}({\rm Ker}\,\rho) = (0)$.
Therefore $m^{\bullet}{\cdot}u = \rho (u)m^{\bullet}$ for all $u
\in U$. The calculation
$$
(m^\perp{\cdot}A)(m) = m^\perp (A{\cdot}m) = m^\perp (km) = (0)
$$
shows that $m^\perp{\cdot}A \subseteq m^\perp$. Therefore
$N^{\bullet} = m^\perp{\cap}M^{\bullet}$ is a right $A$-submodule
of $M^{\bullet}$. Since
$$
(0) \neq m^{\bullet}(M) = m^{\bullet}(H{\cdot}m) =
(m^{\bullet}{\cdot}H)(m) = M^{\bullet}(m)
$$
it follows that $M^{\bullet} \not \subseteq m^\perp$. This means
that the right $A$-submodule $N^{\bullet}$ is a codimension one
subspace of $M^{\bullet}$. Since
$$
(M^{\bullet}{\cdot}({\rm Ker}\,\chi))(m) = M^{\bullet}(({\rm
Ker}\,\chi){\cdot}m) = M^{\bullet}(0) = (0)
$$
we conclude that $M^{\bullet}{\cdot}({\rm Ker}\,\chi) \subseteq
m^\perp{\cap}M^{\bullet} = N^{\bullet}$. We have shown that
 \begin{equation}\label{EqRhoChi}
{\rm ann}_U(km^{\bullet}) = {\rm Ker}\,\rho \qquad \mbox{and}
\qquad {\rm ann}_A(M^{\bullet}/N^{\bullet}) = {\rm Ker}\,\chi.
\end{equation}

We next show that the only right $H$-submodule of $M^{\bullet}$
contained in $N^{\bullet}$ is $(0)$. Let $L$ be a right
$H$-submodule of $M^{\bullet}$ contained in $N^{\bullet}$. Then
$L(M) = L(H{\cdot}m) = (L{\cdot}H)(m) \subseteq m^\perp(m) = (0)$
implies that $L = (0)$. We have shown that $(M^{\bullet},
km^{\bullet}, N^{\bullet}) \in \underline{\mathcal M}'_H$; in
particular $M^{\bullet} \simeq R(\chi, \rho)$.

Consider the bilinear form $\beta : M^{\bullet}{\otimes}M
\longrightarrow k$ given by $\beta (p, n) = p(n)$ for all $p \in
M^{\bullet}$ and $n \in M$. Note that $\beta$ is right
non-singular. By definition of the transpose module action $\beta
(p{\cdot}h, n) = \beta (p, h{\cdot}n)$ for all $p \in
M^{\bullet}$, $h \in H$, and $n \in M$; that is $\beta$ is
$H$-balanced. We observe that $(M^{\bullet})^\perp$ is the largest
$H$-submodule of $M$ contained in $N$. For let $n \in N$. Then
$H{\cdot}n \subseteq N$ if and only if $(0) =
m^{\bullet}(H{\cdot}n) = (m^{\bullet}{\cdot}H)(n) =
M^{\bullet}(n)$; that is $H{\cdot}n \subseteq N$ if and only if
$M^{\bullet}(n) = (0)$. Thus $M/(M^{\bullet})^\perp \simeq L(\rho,
\chi)$ and $\beta$ induces an $H$-balanced bilinear form
$\overline{\beta} : M^{\bullet}{\otimes}(M/(M^{\bullet})^\perp)
\longrightarrow k$. Compare with Proposition \ref{PropBilinear}.
\begin{Prop}\label{PropDuality}
Let $U$ and $A$ be algebras over the field $k$, let $H =
U{\otimes}A$ be an algebra structure on the tensor product of
their underlying vector spaces which satisfies
(\ref{EqUandAInclusions}), and let $(M, km, N) \in
{}_H\underline{\mathcal M}'$.
\begin{enumerate}
\item[{\rm a)}] Let $\chi \in \Alg(A, k)$ and  $\rho \in \Alg(U,
k)$ satisfy
$$
{\rm ann}_A(km) = {\rm Ker}\,\chi \qquad \mbox{and} \qquad {\rm
ann}_U(M/N) = {\rm Ker}\,\rho.
$$
Then $(M^{\bullet}, km^{\bullet}, N^{\bullet}) \in
\underline{\mathcal M}'_H$ and satisfies
$$
{\rm ann}_U(km^{\bullet}) = {\rm Ker}\,\rho \qquad \mbox{and}
\qquad {\rm ann}_A(M^{\bullet}/N^{\bullet}) = {\rm Ker}\,\chi.
$$
\item[{\rm b)}] Suppose that $f : (M, km, N) \longrightarrow (M',
km', N')$ is a morphism in ${}_H\underline{\mathcal M}'$. Then
$f^*(M'^{\bullet}) \subseteq M^{\bullet}$ and the restriction $f^r
= f^*|{M'^{\bullet}}$ is a morphism $f^r : (M^{'\bullet},
m^{'\bullet}, N^{'\bullet}) \longrightarrow (M^{\bullet},
m^{\bullet}, N^{\bullet})$ in $\underline{\mathcal M}'_H$.
\item[{\rm c)}] Suppose that $N^{\bullet}$ is the only
codimension-one right $A$-submodule of $M^{\bullet}$. Then $km$ is
the only one-dimensional left $A$-submodule of $M$.
\end{enumerate}
\end{Prop}

\pf We have shown part a). Part b) is left as an easy exercise. We
establish part c).

Suppose that $n \in M$ and $A{\cdot}n = kn$. Then
$n^\perp{\cap}M^{\bullet} = M^{\bullet}$, in which case
$M^{\bullet}(n) = (0)$, or $n^\perp{\cap}M^{\bullet}$ is a
codimension one right $A$-submodule of $M^{\bullet}$, in which
case $n^\perp{\cap}M^{\bullet} = m^\perp{\cap}M^{\bullet}$ by
assumption. In any event $N^{\bullet}(n) = (0)$.

Let $C = \{ n \in M\, |\, N^{\bullet}(n) = (0)\}$. Observe that $m
\in C$. The rule $f : C \longrightarrow
(M^{\bullet}/N^{\bullet})^*$ given by $f(n)(p + N^{\bullet}) =
p(n)$ is describes a well-defined linear function. We show that
$f$ is one-one. Suppose that $n \in C$ and $f(n) = 0$. Then
$$
(0) = f(n)(M^{\bullet}) = M^{\bullet}(n) =
(m^{\bullet}{\cdot}H)(n) = m^{\bullet}(H{\cdot}n)
$$
implies that $H{\cdot}n \subseteq {\rm Ker}\,m^{\bullet} =
N^{\bullet}$. But $N^{\bullet}$ contains no left $H$-submodules
other than $(0)$. Therefore $n = 0$. We have shown that $f$ is
one-one.

Since $m \in C$, $f$ is one-one, and ${\rm Im}\,f$ is at most
one-dimensional, it follows that $C = km$. This concludes our
proof. \qed
\medskip

The ``right" counterpart of the preceding proposition holds by
virtue of Example \ref{HOp}.
\section{The Main Results for ${}_H\underline{\mathcal M}$}\label{SecMainHM}
We begin by describing the type of algebra of fundamental
importance in \cite{Irr}. These algebras $A$ (and $U$ also) are
generated by a subgroup $\Gamma$ of the group of units of $A$ and
a indexed set of elements $\{a_i\}_{i \in I}$. There is an indexed
set of characters $\{\chi_i\}_{i \in I} \subseteq
\widehat{\Gamma}$ such that $ga_i g^{-1} = \chi_i (g)a_i$ for all
$g \in \Gamma$ and $i \in I$. Let $A'$ be the subalgebra of $A$
generated by $\Gamma$.

Suppose that $\chi_i \neq 1$ for all $i \in I$. Then $\rho (a_i) =
0$ for all $\rho \in \Alg(A, k)$ and $i \in I$. Thus the
restriction map $\Alg(A, k) \longrightarrow \Alg(A', k)$ is
one-one. In important applications $U$ and $A$ below will have
this description and thus the restriction map is one-one.

The theorem of this section is derived from two results.
\begin{Lemma}\label{LemMainU}
Let $U$ and $A$ be algebras over the field $k$ and let $H =
U{\otimes}A$ be an algebra structure on the tensor product of the
underlying vector spaces of $U$ and $A$ which satisfies
(\ref{EqUandAInclusions}). Suppose that $U'$ is a subalgebra of
$U$ such that:
\begin{enumerate}
\item[{\rm a)}] The restriction map $\Alg(U, k) \longrightarrow
\Alg(U', k)$ is one-one \item[{\rm b)}] and
$(u{\otimes}a){\cdot}_\chi u' = \chi (a)uu'$ for all $u \in U$, $a
\in A$, $\chi \in \Alg(A, k)$, and $u' \in U'$.
\end{enumerate}
Let $\rho, \rho' \in \Alg(U, k)$ and $\chi, \chi' \in \Alg(A, k)$.
If $L(\rho, \chi) \simeq L(\rho', \chi')$ as left $U$-modules then
$\rho = \rho'$. In particular there is a unique codimension one
left $U$-module of  $L(\rho, \chi)$ which contains $I(\rho,
\chi)$.
\end{Lemma}

\pf Suppose that $L(\rho, \chi) \simeq L(\rho', \chi')$ as left
$U$-modules and consider the composite of left $U$-modules $f : U
\longrightarrow L(\rho, \chi) \simeq L(\rho', \chi') = M$. Since
the latter contains a codimension one left $U$-submodule $N$ with
${\rm ann}_U(M/N) = {\rm Ker}\,\rho'$, it follows that $I(\rho,
\chi) \subseteq f^{-1}(N) = {\rm Ker}\,\rho'$. We have shown that
$I(\rho, \chi) \subseteq {\rm Ker}\rho'$. The calculation
$\rho((u{\otimes}a){\cdot}_\chi u') = \rho(\chi (a)uu') = \chi
(a)\rho(u)\rho(u')$ for all $u \in U$, $a \in A$, and $u' \in U'$
shows that $H{\cdot}_\chi ({\rm Ker}\,\rho{\cap}U') \subseteq {\rm
Ker} \rho$. Therefore
$$
({\rm Ker}\,\rho){\cap}U' \subseteq I(\rho, \chi) \subseteq {\rm
Ker} \rho'
$$
which implies that $({\rm Ker}\,\rho){\cap}U' = ({\rm
Ker}\,\rho'){\cap}U'$ since both intersections are codimension one
subspaces of $U'$. The preceding equation implies $\rho|{U'} =
\rho'|{U'}$ from which $\rho = \rho'$ follows by assumption. The
last statement in the conclusion of the lemma is evident. \qed
\begin{Lemma}\label{LemMainA}
Let $U$ and $A$ be algebras over the field $k$ and let $H =
U{\otimes}A$ be an algebra structure on the tensor product of the
underlying vector spaces of $U$ and $A$ which satisfies
(\ref{EqUandAInclusions}). Suppose that $A'$ is a subalgebra of
$A$ such that:
\begin{enumerate}
\item[{\rm a)}] The restriction map $\Alg(A, k) \longrightarrow
\Alg(A', k)$ is one-one and \item[{\rm b)}] $a'{\cdot}_\rho
(u{\otimes}a) = \rho (u)a'a$ for all $a' \in A'$, $\rho \in
\Alg(U, k)$, $u \in U$, and $a \in A$.
\end{enumerate}
Let $\rho \in \Alg(U, k)$ and $\chi, \chi' \in \Alg \,(A, k)$. If
$L(\rho, \chi) \simeq L(\rho, \chi')$ as left $H$-modules then
$\chi = \chi'$.
\end{Lemma}

\pf Regard $k$ as a left $U$-module via $u{\cdot}1 = \rho(u)$ for
all $u \in U$. Consider the composite $\rho' : U_\chi
\longrightarrow k$ of left $U$-module maps given by $U_\chi
\longrightarrow L(\rho, \chi) \stackrel{f}{\longrightarrow}
L(\rho, \chi') \stackrel{\overline{\rho}}{\longrightarrow} k$,
where the first map is the projection, $f$ is an isomorphism of
left $H$-modules, and the third $\overline{\rho}$ is the map is
given by $u + I(\rho, \chi') \mapsto \rho(u)$ for all $u \in U$.
Since $\rho'$ is a left $U$-module map we have $\rho'(u) =
\rho'(u1) = u{\cdot}\rho'(1) = \rho(u)\rho'(1)$. Therefore $\rho'
= \rho'(1)\rho$.

Let $u \in U$ satisfy $f(1 + I(\rho, \chi)) = u + I(\rho, \chi')$
and let $a \in A$. Using the definition of $\rho'$, the fact that
$f$ is a map of left $A$-modules, and part a) of Proposition
\ref{PropBilinear}, we see that
\begin{eqnarray*}
\rho'((1{\otimes}a){\cdot}_\chi (1 + I(\rho, \chi)) & = & \rho
(f((1{\otimes}a)
{\cdot}_\chi (1 + I(\rho, \chi))) \\
& = & \overline{\rho}((1{\otimes}a){\cdot}_{\chi'}f(1 + I(\rho, \chi))) \\
& = & \overline{\rho} ((1{\otimes}a){\cdot}_{\chi'}(u + I(\rho, \chi'))) \\
& = & \rho ((1{\otimes}a){\cdot}_{\chi'}u) \\
& = & \chi' (a{\cdot}_\rho (u{\otimes}1)).
\end{eqnarray*}

Since $\rho' = \rho'(1)\rho$ we calculate on the other hand that
\begin{eqnarray*}
\rho'((1{\otimes}a){\cdot}_\chi (1 + I(\rho, \chi))
& = & \rho'(1) \rho((1{\otimes}a){\cdot}_\chi (1 + I(\rho, \chi))  \\
& = & \rho'(1) \rho((1{\otimes}a){\cdot}_\chi 1) \\
& = & \rho'(1)\chi (a).
\end{eqnarray*}
We have shown that $\chi'(a{\cdot}_\rho (u{\otimes}1)) =
\rho'(1)\chi (a)$ for all $a \in A$. Now suppose that $a' \in A'$.
By virtue of the preceding equation
$$
\chi'(a')\rho(u)  = \chi'(a'\rho(u)\chi'(1)) =
\chi'(a'{\cdot}_\rho (u{\otimes}1)) = \rho'(1)\chi (a');
$$
the second equation  follows by assumption. Therefore $\rho'(u) =
1$ and $\chi(a') = \chi'(a')$ for all $a' \in A'$. By assumption
$\chi = \chi'$. \qed
\medskip

For a category $\mathcal{M}$ we denote the isomorphism classes of
objects in $\mathcal{M}$ by $[\mathcal{M}]$ and for an object $M
\in \mathcal{M}$ we let $[M]$ be the isomorphism class of $M$. As
a consequence of the two preceding lemmas:
\begin{Theorem}\label{ThmMain}
Let $U$ and $A$ be algebras over the field $k$ and let $H =
U{\otimes}A$ be an algebra structure on the tensor product of the
underlying vector spaces of $U$ and $A$ which satisfies
(\ref{EqUandAInclusions}). Suppose that:
\begin{enumerate}
\item[{\rm a)}] $U'$ is a subalgebra of $U$,
$(u{\otimes}a){\cdot}_\chi u' = \chi (a)uu'$ for all $u \in U$, $a
\in A$, $\chi \in \Alg(A, k)$, and $u' \in U'$, and the
restriction map $\Alg(U, k) \longrightarrow \Alg(U', k)$ is
one-one; \item[{\rm b)}] $A'$ is a subalgebra of $A$,
$a'{\cdot}_\rho (u{\otimes}a) = \rho (u)a'a$ for all $a' \in A'$,
$\rho \in \Alg(U, k)$, $u \in U$, and $a \in A$, and the
restriction map $\Alg(A, k) \longrightarrow \Alg(A', k)$ is
one-one.
\end{enumerate}
Then $\Alg(U, k){\times}\Alg(A, k) \longrightarrow
[{}_H\underline{\mathcal M}]$ given by $(\rho, \chi) \mapsto
[L(\rho, \chi)]$ is bijective. \qed
\end{Theorem}
\medskip

Observe that the hypothesis of the theorem holds for
$H^{\tilde{op}}$ as well, where $A'^{op}$ in $H^{\tilde{op}}$
plays the role of $U'$ in $H$ and $U'^{op}$ in $H^{\tilde{op}}$
plays the role of $A'$ in $H$. Therefore:
\begin{Cor}\label{CorMainRight}
Under the hypothesis of the preceding theorem, the function
$\Alg(A, k){\times}\Alg(U, k) \longrightarrow [\underline{\mathcal
M}_H]$ given by $(\chi, \rho) \mapsto [R(\chi, \rho)]$ is
bijective. \qed
\end{Cor}
The theorem has interesting consequences for objects of the
category ${}_H\underline{\mathcal M}$.
\begin{Cor}\label{CorUniqueSubModules}
Let $U$ and $A$ be algebras over the field $k$ and let $H =
U{\otimes}A$ be an algebra structure on the tensor product of the
underlying vector spaces of $U$ and $A$ which satisfies
(\ref{EqUandAInclusions}). Assume that the hypothesis of Theorem
\ref{ThmMain} holds. Then:
\begin{enumerate}
\item[{\rm a)}] Every object of ${}_H\underline{\mathcal M}$ has a
unique one-dimensional $A$-submodule and a unique codimension one
$U$-submodule. \item[{\rm b)}] Suppose that $f : M \longrightarrow
M'$ is a left $H$-module map, where $M$ and $M'$ are objects of
${}_H\underline{\mathcal M}$. Then $f = 0$ or $f$ is an
isomorphism.
\end{enumerate}
\end{Cor}

\pf We first show part a). We have noted that the hypothesis of
the theorem apples to $H^{\tilde{op}}$. In light of Lemma
\ref{LemMainU} and Theorem \ref{ThmMain} we need only show that
the object $M \in \underline{\mathcal M}$ has a $1$-dimensional
$A$-submodule. Let $M \in {}_H\underline{\mathcal M}$ and let $(M,
km, N) \in {}_H\underline{\mathcal M}$ be derived from $M$. Let
$\rho \in \Alg(U, k)$ and $\chi \in \Alg(A, k)$ satisfy
(\ref{EqChiRho}). Regard the right $H$-module $M^{\bullet}$ as a
left $H^{\tilde{op}}$-module by pullback along the algebra map
$H^{\tilde{op}} \longrightarrow H^{op}$ given by $a{\otimes}u
\mapsto u{\otimes}a$ for all $a \in A$ and $u \in U$. Then
$M^{\bullet} \simeq L^{\tilde{op}}(\chi, \rho)$ as left
$H^{\tilde{op}}$-modules, where $L^{\tilde{op}}(\chi, \rho)$  is
the counterpart of $L(\rho, \chi)$ for $H$. By Lemma
\ref{LemMainU} and Theorem \ref{ThmMain} there is only one
codimension left $H^{\tilde{op}}$-module, or equivalently right
$H$-module, in $M^{\bullet}$. Therefore $M$ has a unique
one-dimensional left $A$-module by part c) of Proposition
\ref{PropDuality}.

Part b) follows from part a). We first note that since $f$ is a
map of left $H$-modules it is also a map of left $A$-modules and
left $U$-modules.

Let $km$ be a one-dimensional left $A$-submodule of $M$ which
generates $M$ as a left $H$-module. Since $f$ is a map of left
$H$-modules $f(M) = f(H{\cdot}m) = H{\cdot}f(m)$. Since $f$ is a
map of left $A$-modules $kf(m)$ is a left $A$-submodule of $M'$.

Suppose that $f \neq 0$. Then $f(m) \neq 0$. Therefore $kf(m)$ is
a one-dimensional left $A$-submodule of $M'$. Now $M'$ is
generated as a left $H$-module by some one-dimensional left
$A$-submodule of $M'$. By uniqueness this submodule must be
$kf(m)$. Thus $f$ is onto. It remains to show that $f$ is one-one.

Now $M'$ contains a codimension one left $U$-submodule $N'$. Since
$f$ is onto and a map of left $U$-modules $f^{-1}(N')$ is a
codimension one left $U$-submodule of $M$. Now $M$ has a
codimension one left $U$-submodule which contains no left
$H$-submodule of $M$ other than $(0)$. By uniqueness this
submodule must be $f^{-1}(N')$. Since ${\rm Ker}\,f \subseteq
f^{-1}(N')$ and is an $H$-submodule of $M$ it follows that ${\rm
Ker}\,f = 0$. We have shown that $f$ is one-one. \qed
\medskip

For an algebra $A$ we denote the  set of isomorphism classes of
finite-dimensional irreducible left $A$-modules by $\Irr(A)$. As a
result of the preceding corollary:
\begin{Cor}\label{CorSimples}
Let $U$ and $A$ be algebras over the field $k$ and let $H =
U{\otimes}A$ be an algebra structure on the tensor product of the
underlying vector spaces of $U$ and $A$ which satisfies
(\ref{EqUandAInclusions}). Assume that the hypothesis of Theorem
\ref{ThmMain} holds and also that the irreducible left $U$-modules
and irreducible left $A$-modules are one-dimensional. Then:
\begin{enumerate}
\item[{\rm a)}] The finite-dimensional irreducible left
$H$-modules are the same as the finite-dimensional objects of
${}_H\underline{\mathcal M}$. \item[{\rm b)}] Suppose that $U$ and
$A$ are finite-dimensional. Then the function $\Alg(U,
k){\times}\Alg(A, k) \to \Irr(H)$ given by $(\chi, \rho) \mapsto
[L(\rho, \chi)]$ is bijective.
\end{enumerate}
\qed
\end{Cor}
\section{When $U$ and $A$ are Bialgebras}\label{SecMainBialg}
Let $U$ and $A$ be bialgebras over the field $k$, suppose $\tau :
U{\otimes}A \longrightarrow k$ is convolution invertible and
satisfies (A.1)--(A.4), and let $H = (U{\otimes}A)^{\sigma}$. In
this section we apply the major ideas of the preceding section to
$H$.

Suppose that $\rho \in G(U^o) = \Alg(U, k)$ and $\chi \in G(A^o) =
\Alg(A, k)$. Using (\ref{Eq2CocyleMult}) we see that
(\ref{EqUChiAction}) in this case is
\begin{equation}\label{EqbialgUHChi}
(u{\otimes}a){\cdot}_\chi u' = u\tau(u'_{(1)},
a_{(1)})u'_{(2)}\chi(a_{(2)})\tau^{-1}(u'_{(3)}, a_{(3)})
\end{equation}
for all $u, u' \in U$ and $a \in A$ and (\ref{EqARhoAction}) in
this case is
\begin{equation}\label{EqbialgUHRho}
a{\cdot}_\rho(u{\otimes}a') = \tau(u_{(1)},
a_{(1)})\rho(u_{(2)})a_{(2)}\tau^{-1}(u_{(3)}, a_{(3)})a'
\end{equation}
for all $a, a' \in A$ and $u \in U$.  Observe that
\begin{equation}\label{EqPhiUABialg}
\Psi = (\tau(\rho{\otimes}\chi)\tau^{-1}){\circ}\tau_{A, U}
\end{equation}
is just conjugation of $\rho{\otimes}\chi$ by $\tau$ in the dual
algebra $(U{\otimes}A)^*$ preceded by the twist map.

Now let $U' = kG(U)$, $A' = kG(A)$, and suppose that $G(U^o)$ and
$G(A^o)$ are commutative groups. Using (\ref{EqbialgUHChi}) we see
for $u \in U$, $a \in A$, $\chi \in \Alg(A, k)$, and $u' \in G(U)$
that
\begin{eqnarray*}
(u{\otimes}a){\cdot}_\chi u'  & = & u\tau (u',
a_{(1)})u'\chi (a_{(2)})\tau^{-1}(u', a_{(3)}) \\
& = & u\left(\PH (\tau_\ell(u') \chi
\tau_\ell(u')^{-1})(a)\right)u'  \\ & = & uu'\chi (a)
\end{eqnarray*}
and therefore $(u{\otimes}a){\cdot}_\chi u' = uu'\chi (a)$ for all
$u \in U, a \in A$ and $u' \in U'$. Likewise using
(\ref{EqbialgUHRho}) it follows that $a'{\cdot}_\rho (u{\otimes}a)
= a'a\rho(u)$ for all $a' \in A'$, $\rho \in \Alg(U, k)$, $u \in
U$, and $a \in A$.
\begin{Theorem}\label{MainBialg}
Let $U$ and $A$ be bialgebras over the field $k$, suppose $\tau :
U{\otimes}A \longrightarrow k$ is convolution invertible and
satisfies (A.1)--(A.4), and let $H = (U{\otimes}A)^{\sigma}$.
Suppose that all $\rho \in \Alg(U,k)$, $\chi \in \Alg(A,k)$ are
determined by their respective restrictions $\rho|{G(U)}$,
$\chi|{G(A)}$. Then $\Alg(U,k)$, $\Alg(A,k)$ are abelian groups
and:
\begin{enumerate}
\item[{\rm a)}] $\Alg(U,k) {\times}\Alg(A,k)\longrightarrow
[{}_H\underline{\mathcal M}]$ given by $(\rho, \chi) \mapsto
[L(\rho, \chi)] $ is bijective. \item[{\rm b)}] Every object of
${}_H\underline{\mathcal M}$ has a unique one-dimensional
$A$-submodule and a unique codimension one $U$-submodule.
\item[{\rm c)}] Suppose that the finite-dimensional irreducible
left $U$-modules and the finite-dimensional irreducible left
$A$-modules are one-dimensional. Then the finite-dimensional
irreducible left $H$-modules are the finite-dimensional objects of
${}_H\underline{\mathcal M}$.
\end{enumerate}
\end{Theorem}

\pf The hypothesis of Theorem \ref{ThmMain} holds for $H$ with $U'
= kG(U)$ and $A' = kG(A)$. Thus part a) follows. Part b) is part
a) of Corollary \ref{CorUniqueSubModules} and part c) is part a)
of Corollary \ref{CorSimples}. \qed
\medskip

We will show that the preceding theorem applies to a wide class of
pointed Hopf algebras. First we recall a basic Hopf module
construction.

Let $A$ be a Hopf algebra with sub-Hopf algebra $B$. Suppose that
$D, C$ are subcoalgebras of $A$ which satisfy $D \subseteq C$, $BD
\subseteq D$, $BC \subseteq C$, and $\Delta (C) \subseteq
B{\otimes}C + C{\otimes}D$. Set $M = C/D$ and write $\overline{c}
= c + C$ for all $c \in C$. Then $M$ is a left $B$-Hopf module,
where $b{\cdot}\overline{c} = \overline{bc}$ and $\rho
(\overline{c}) = c_{(1)}{\otimes}\overline{c_{(2)}}$ for all $c
\in C$. By the Fundamental Theorem for Hopf modules \cite[Theorem
4.1.1]{SweedlerBook} it follows that $M$ is $(0)$ or a free left
$B$-module with basis any linear basis of $M^{\co B} =\{
\overline{c} \, | \, \rho (\overline{c}) =
1{\otimes}\overline{c}\}$.

Now suppose that $B = A_0$ is a sub-Hopf algebra of $A$ and let $n
> 0$. Then $C = A_n$ and $D = A_{n-1}$ satisfy the conditions of the
preceding paragraph. Therefore $A_n/A_{n-1}$ is $(0)$ or a free
left $B$-module. As a consequence $A_n$ is a free left $B$-module
for all $n \geq 0$; thus $A$ is a free left $B$-module.

Any two bases for a free module over a Hopf algebra $B$ have the
same cardinality, and thus rank of the free module is
well-defined, since Hopf algebras are augmented algebras. For the
same reason, if $M$ is a free left $B$-submodule of a free left
$B$-module $N$ then ${\rm rank}\, M \leq {\rm rank}\,N$.

Suppose that $a$ is a skew primitive element of $A$. We say that
$a$ is of {\em finite type} if the sub-Hopf algebra of $A$
generated by $A_0{\cup}\{a\}$ is a free left $A_0$-module of
finite rank.
\begin{Cor}\label{CorSkewPrimFR}
Let $U$ and $A$ be pointed Hopf algebras over and algebraically
closed field $k$ of characteristic zero and suppose $\tau :
U{\otimes}A \longrightarrow k$ is convolution invertible and
satisfies (A.1)--(A.4). Suppose further that $U$ and $A$ are
generated by skew primitives of finite rank and have commutative
coradicals. Then the conclusions of the preceding theorem hold for
$H = (U{\otimes}A)^{\sigma}$.
\end{Cor}

\pf We need only show that all $\rho \in G(U^o)$, $\chi \in
G(A^o)$ are determined by their respective restrictions
$\rho|{G(U)}$, $\chi|{G(A)}$. We give an argument for $A$ which is
automatically an argument for $U$ also. To this end we need only
show that there is an indexed set of skew primitive elements
$\{a_i\}_{i \in I}$, which together with $\Gamma$ generate $A$ as
an algebra, and there is an indexed set of non-trivial characters
$\{\chi\}_{i \in I}$ such that $ha_i h^{-1} = \chi_i(h)a_i$ for
all $h \in \Gamma$ and $i \in I$. See the opening commentary for
Section \ref{SecMainHM}.

Let $B = A_0$. Let $\Gamma = G(A)$. Since $A_0$ is cocommutative
and $k$ is algebraically closed it follows that $B = k\G$. Since
$B$ is commutative $\Gamma$ is a commutative group.

Suppose that $a \in A$ is a skew primitive element of finite rank.
We may assume that $\Delta (a) = a{\otimes}g + 1{\otimes}a$ for
some $g \in \Gamma$ and that $a \not \in B$. Let $E$ be the
sub-Hopf algebra of $A$ generated by $B{\cup}\{a\}$. Then $E$ is a
free left $B$-module of finite rank by assumption.

Let $V = \{ v \in E \, | \, \Delta (v) = v{\otimes}g +
1{\otimes}v\}$. Then $C = BVB$ is a left $B$-module, a
subcoalgebra of $E$, and $\Delta (C) \subseteq C{\otimes}B +
B{\otimes}C$. Thus $M = C/B$ is a left $B$-Hopf module. Since $1 +
{\rm rank}\,M = {\rm rank}\,C \leq {\rm rank}\,E$, and the latter
is finite, it follows that ${\rm rank}\,M$ is finite. Now
$\overline{V} \subseteq M^{\co B}$. Thus ${\rm Dim}\,\overline{V}
\leq {\rm rank}\,M$ is finite. Since $V{\cap}B = k(g-1)$ we
conclude that $V$ is a finite-dimensional vector space.

Let $h \in \Gamma$. Since $\Gamma$ is commutative and $E$ is a
subalgebra of $A$ which contains $h$ it follows that $hVh^{-1}
\subseteq V$. By assumption $V \neq (0)$. Since $V$ is
finite-dimensional and $k$ is an algebraically closed field of
characteristic zero, there is a basis $v_1, \ldots, v_n$ for $V$
consisting of common eigenvectors for the conjugation action by
$\Gamma$. Therefore there are characters $\chi_1, \ldots, \chi_n
\in \widehat{\Gamma}$ such that $hv_i h^{-1} = \chi_i (h)v_i$ for
all $h \in \Gamma$ and $1 \leq i \leq n$.

Fix $1 \leq i \leq n$. Observe that $v_i$ and $\Gamma$ generate
sub-Hopf algebra of $A$ of finite rank. By calculations found in
\cite[Section 3]{Simpoint} it follows that $\chi_{i} (g) = 1$
implies that $v_i \in B$. Since $a$ is in the span of the $v_i$'s,
it is clear how to form the families $\{a_i\}_{i \in I}$ and
$\{\chi_i\}_{i \in I}$ which satisfy the conditions outlined at
the beginning of the proof. \qed
\medskip

For later use we note:
\begin{Lemma}\label{modulo}
Let $U$ and $A$ be bialgebras over the field $k$ and assume that
$G(U^0)$ is abelian. Suppose $\tau : U{\otimes}A \longrightarrow
k$ is convolution invertible and satisfies (A.1)--(A.4), and let
$H = (U{\otimes}A)^{\sigma}$. Let $u\in U$ and $g \in G(A)$ and
assume that $u \otimes g$ is central in $H.$ Then for all $\rho
\in \Alg(U,k)$ and $\chi \in \Alg(A,k)$ the following are
equivalent:
\begin{enumerate}
\item [{\rm a)}] $u \otimes g - 1 \otimes 1$ acts as zero on
$L(\rho,\chi).$ \item [{\rm b)}] $\rho(u)\chi(g) = 1.$
\end{enumerate}
\end{Lemma}
\pf Part a) implies
\begin{equation}\label{a)}
(u \otimes g - 1 \otimes 1)\cdot_{\chi}v = u v\sw2\tau(v\sw1,g)
\chi(g) \tau^{-1}(v\sw3,g) \in I(\rho,\chi)
\end{equation}
for all $v \in U$. If $v=1$ then \eqref{a)} implies $u\chi(g)- 1
\in I(\rho,\chi),$ hence $\rho(u\chi(g)- 1)=0$ or equivalently
$\rho(u)\chi(g) = 1.$ Thus part a) implies part b).

Conversely, assume part b). Since $u \otimes g$ is central in $H$,
the $k$-span of $(u \otimes g - 1 \otimes 1)\cdot_{\chi}v, v \in
U,$ is an $H$-submodule of $U_{\chi}.$ Moreover, for any $v \in
U,$
$$\tau(v\sw1,g)\rho(v\sw2) \tau^{-1}(v\sw3,g) = \rho(v),$$
since $\tau_r(g) \in \Alg(U,k)$ and $G(U^0)=\Alg(U,k)$ is abelian.
Hence
\begin{align*}
\rho((u \otimes g)\cdot_{\chi}v)&= \rho(u v\sw2 \tau(v\sw1,g) \chi(g) \tau^{-1}(v\sw3,g))\\
 &= \rho(u) \chi(g) \rho(v)\\
 &= \rho(v)
 \end{align*}
 by b). This proves \eqref{a)} by definition of $I(\rho,\chi).$
\epf
\medskip

In connection with Propositions \ref{PropfLift} and
\ref{PropgLift} we will be interested in bialgebra maps of
bialgebras of the type $(U{\otimes}A)^\sigma$.
\begin{Prop}\label{PropUABialgebraMorp}
Let $U$, $\overline{U}$, $A$, and $\overline{A}$ be bialgebras
over $k$, suppose $\tau : U{\otimes}A \longrightarrow k$ and
$\overline{\tau} : \overline{U}{\otimes}\overline{A}
\longrightarrow k$ satisfy (A.1)--(A.4), and suppose that $f : U
\longrightarrow \overline{U}$ and $g : A \longrightarrow
\overline{A}$ are bialgebra maps which satisfy
$\overline{\tau}{\circ}(f{\otimes}g) = \tau$. Then $f{\otimes}g :
(U{\otimes}A)^\sigma \longrightarrow
(\overline{U}{\otimes}\overline{A})^{\overline{\sigma}}$ is a
bialgebra map.
\end{Prop}

\pf Since $f$ and $g$ are coalgebra maps $f{\otimes}g :
U{\otimes}A \longrightarrow \overline{U}{\otimes}\overline{A}$ is
a coalgebra map of the tensor product of coalgebras. As the
underlying coalgebra structures of $(U{\otimes}A)^\sigma$ and
$(\overline{U}{\otimes}\overline{A})^{\overline{\sigma}}$ are
$U{\otimes}A$ and $\overline{U}{\otimes}\overline{A}$
respectively, it follows that $f{\otimes}g : (U{\otimes}A)^\sigma
\longrightarrow
(\overline{U}{\otimes}\overline{A})^{\overline{\sigma}}$ is a
coalgebra map. Since $(f{\otimes}g)^* :
(\overline{U}{\otimes}\overline{A})^* \longrightarrow
(U{\otimes}A)^*$ is an algebra map $\overline{\tau}^{-1} =
((f{\otimes}g)^*(\tau))^{-1} = (f{\otimes}g)^*(\tau^{-1}) =
\tau^{-1}{\circ}(f{\otimes}g)$. At this point is it easy to see
that $f{\otimes}g : (U{\otimes}A)^\sigma \longrightarrow
(\overline{U}{\otimes}\overline{A})^{\overline{\sigma}}$ is an
algebra map. \qed

Suppose that $U$ and $A$ are bialgebras over $k$ and $\tau :
U{\otimes}A \longrightarrow k$ is a linear form which satisfies
(A.1)--(A.4). Then $\tau$ is convolution invertible by Lemma
\ref{LemmaTauInver} if $U$ or $A^{op}$ is a Hopf algebra, in
particular if $A$ is a Hopf algebra with bijective antipode.

Suppose in addition that $A$ has bijective antipode, set $H =
(U{\otimes}A)^\sigma$, and let $f : H \longrightarrow D(A)$ be the
bialgebra map defined at the end of Section \ref{SecPrelim} by
$f(u{\otimes}a) = \tau_\ell(u){\otimes}a$ for all $u \in U$ and $a
\in A$. Let $\rho \in G(U^o)$ and $\chi \in G(A^o)$. We will
examine $L(\rho, \chi)$ in the context of part c) of Proposition
\ref{PropBilinear} and find a condition for there to be a left
$D(A)$-module such that pullback along $f$ explains $L(\rho,
\chi)$. To understand modules for $D(A)$ we need to review a
variant of the Yetter--Drinfeld category discussed in many places;
in particular in \cite{ASSurvey}.

Let $B$ be any bialgebra over $k$ and let ${}_B{\mathcal YD}^B$ be
the category whose objects are triples $(M, {\cdot}, \delta)$,
where $(M, {\cdot})$ is a left $B$-module and $(M, \delta)$ is a
right $B$-comodule which are compatible in the sense
\begin{equation}\label{EqYDLeftRight}
b_{(1)}{\cdot}m_{(0)}{\otimes}b_{(2)}m_{(1)} =
(b_{(2)}{\cdot}m)_{(0)}{\otimes} (b_{(2)}{\cdot}m)_{(1)}b_{(1)}
\end{equation}
for all $b \in B$ and $m \in M$, and whose morphisms are maps of
left $B$-modules and right $B$-comodules. We observe that when
$B^{op}$ has antipode $T$ then (\ref{EqYDLeftRight}) is equivalent
to
\begin{equation}\label{EqYDLeftRightVarsigma}
\delta (b{\cdot}m) =
b_{(2)}{\cdot}m_{(0)}{\otimes}b_{(3)}m_{(1)}T(b_{(1)})
\end{equation}
for all $b \in B$ and $m \in M$. An example, which is the
centerpiece of \cite{RadRec} in the study of simple modules for
the double, is the following \cite[Lemma 2]{RadRec}:
\begin{Ex}\label{ExSimpleYD}
Let $B$ be a bialgebra over $k$, suppose that $B^{op}$ is a Hopf
algebra with antipode $T$, and let $\beta \in G(B^o)$. Then $(B,
{\succ}_\beta, \Delta) \in {}_B{\mathcal YD}^B$, where
$$
b{\succ}_\beta m = (b_{(2)}{\leftharpoonup}\beta)mT(b\sw1)
$$
for all $b, m \in B$.
\end{Ex}

The map $b \mapsto b{\leftharpoonup}\beta$ of the example is an
algebra automorphism of $B$. Thus the module $(B, {\succ}_\beta)$
can be regarded as a generalized adjoint action.

Now suppose that $B$ is a Hopf algebra with bijective antipode $S$
and let $(M, {\cdot}, \delta) \in {}_B{\mathcal YD}^B$. Then $(M,
{\bullet}) \in {}_{D(B)}{\mathcal M}$ where
$$
(p{\otimes}b){\bullet}m = p{\rightharpoonup}(b{\cdot}m) =
(b{\cdot}m)_{(0)}{<}p, (b{\cdot}m)_{(1)}{>}
$$
for all $p \in B^o$, $b \in B$, and $m \in M$. When $B$ is
finite-dimensional the preceding equation describes the essence of
a categorical isomorphism of ${}_{D(B)}{\mathcal M}$ and
${}_B{\mathcal YD}^B$. See the primary reference \cite{Majid} as
well as \cite[Proposition 3.5.1]{AAQYBE}.

It is also interesting to note that ${}_{B^o}{\mathcal YD}^{B^o}$
also accounts for left modules for $D(B)$. For suppose that $(M,
{\cdot}, \delta) \in {}_{B^o}{\mathcal YD}^{B^o}$ and let $i_B : B
\longrightarrow (B^o)^*$ be the algebra map defined by $i_B(b)(p)
= p(b)$ for all $b \in B$ and $p \in B^o$. Set $i = i_B$. Then
$(M, {\bullet}) \in {}_{D(B)}{\mathcal M}$, where
$$
(p{\otimes}b){\bullet}m = p{\cdot}(i(b){\rightharpoonup}m) =
p{\cdot}m_{(0)}{<}m_{(1)}, b{>}
$$
for all $p \in B^o$, $b \in B$, and $m \in M$. Observe that the
action on $(M, {\bullet})$ restricted to $B$ is locally finite.
The preceding equation describes the essence of a categorical
isomorphism between the full subcategory of ${}_{D(B)}{\mathcal
M}$ whose objects are locally finite as left $B$-modules and
${}_{B^o}{\mathcal YD}^{B^o}$. Thus when $B$ is finite-dimensional
there is a categorical isomorphism of ${}_{D(B)}{\mathcal M}$ and
${}_{B^o}{\mathcal YD}^{B^o}$. It is the Yetter--Drinfeld category
${}_{B^o}{\mathcal YD}^{B^o}$ which is most appropriate here.

Using Lemma \ref{LemmaTauInver} and (\ref{EqPhiUABialg}) we have
\begin{equation}\label{EqPhiEllUHHopf}
\Psi_r(u) = \tau_\ell (\rho{\rightharpoonup}u_{(1)})\chi
S^{-1}(\tau_\ell(u_{(2)})) = \tau_\ell (u_{(1)})\chi
S^{-1}(\tau_\ell(u_{(2)}){\leftharpoonup}\rho)
\end{equation}
for all $u \in U$.

Regard $A^*$ as a left $H$-module with the transpose action
arising from the right $H$-module structure $(A, {\cdot}_\rho)$.
By part c) of Proposition \ref{PropBilinear} there is an
isomorphism of left $H$-modules $F : L(\rho, \chi) \longrightarrow
U{\cdot}_\rho\chi \subseteq A^*$ given by $F(\overline{u}) =
\Psi_r(u) = u{\cdot}_\rho\chi$ for all $u \in U$, where
$\overline{u} = u + I(\rho, \chi)$. Thus
\begin{equation}\label{EqLRhoChiHopf}
F(\overline{u}) = \tau_\ell (\rho{\rightharpoonup}u_{(1)})\chi
S^{-1}(\tau_\ell(u_{(2)})) = u{\cdot}_\rho\chi
\end{equation}
for all $u \in U$ by (\ref{EqPhiUABialg}). In particular
$U{\cdot}_\rho\chi = {\rm Im}\,F \subseteq A^o$.

The $H$-module $L(\rho, \chi)$ can be explained in terms of $D(A)$
when a very natural condition is satisfied. Let $i = i_A$. Then $i
(g) \in \Alg(A^o, k) = G((A^o)^o)$ and the calculation
\begin{eqnarray*}
a{\cdot}_\rho(u{\otimes}a')
& = & \tau(u_{(1)}, a_{(1)})\rho(u_{(2)})a_{(2)}\tau^{-1}(u_{(3)}, a_{(3)})a' \\
& = & \tau(u_{(1)}, a_{(1)})\tau(u_{(2)}, g)a_{(2)}\tau^{-1}(u_{(3)}, a_{(3)})a' \\
& = & \tau(u_{(1)(1)}, a_{(1)})\tau(u_{(1)(2)}, g)a_{(2)}\tau(u_{(2)}, s^{-1}(a_{(3)}))a' \\
& = &
\left(\PH\tau_\ell(u_{(1)})_{(2)}(a_{(1)})\right)\left(\PH\tau_\ell(u_{(1)})_{(1)}(g)
\right) a_{(2)}\left(\PH S^{-1}(\tau_\ell(u_{(2)}))(a_{(3)})\right)a' \\
& = & \left(\PH
(\tau_\ell(u_{(1)}){\leftharpoonup}i(g))(a_{(1)})\right)a_{(2)}\left(\PH
S^{-1}(\tau_\ell(u_{(2)}))(a_{(3)})\right)a'
\end{eqnarray*}
shows that
\begin{equation}\label{EqARhoRightDbl}
a{\cdot}_\rho(u{\otimes}a')  = \left(\PH
S^{-1}(\tau_\ell(u)_{(1)}){\rightharpoonup}a{\leftharpoonup}(\tau_\ell(u)_{(2)}
{\leftharpoonup}i(g))\right)a'
\end{equation}
for all $a, a\in A$ and $u \in U$. As a consequence
\begin{equation}\label{EqARhoRightDblDual}
(u{\otimes}a'){\cdot}_\rho p = \left(\PH
\tau_\ell(u)_{(2)}{\leftharpoonup}i(g)\right)\left( \PH
i(a'){\rightharpoonup}p\right)\left(\PH S^{-1}(\tau_\ell
(u)_{(1)})\right)
\end{equation}
for all $u \in U$, $a' \in A$, and $p \in A^o$.

Let $(A^o, {\succ}_{i(g)}, \Delta)$ be the object of
${}_{A^o}{\mathcal YD}^{A^o}$ defined in Example \ref{ExSimpleYD}
and let $(M, {\bullet})$ be associated left $D(A)$-module
structure. Using (\ref{EqARhoRightDblDual}) it is not hard to see
that
$$
(u{\otimes}a'){\cdot}_\rho p = f(u{\otimes}a'){\bullet}p
$$
for all $u \in U$, $a' \in A$, and $p \in A^o$. Since $\chi \in
G(A^o)$ we conclude that $H{\cdot}_\rho\chi = U{\cdot}_\rho\chi
\subseteq A^o{\bullet}\chi$, the latter is a $D(A)$-submodule of
$A^o$, and that the map $F : L(\rho, \chi) \longrightarrow
A^o{\bullet}\chi$ given by $\overline{u} \mapsto \Psi_r(u) =
u{\cdot}_\rho \chi$ is a one-one map of left $H$-modules, where
$A^o{\bullet}\chi$ has the left $H$-module structure action
obtained by pullback along $f$.

Note that
\begin{equation}\label{EqPhiRUHHopf}
\Psi_\ell(a) = \tau_r(a_{(1)})\rho
S(\tau_r(a_{(2)}{\leftharpoonup}\chi))
\end{equation}
follows for all $a \in A$ by (\ref{EqPhiUABialg}) also, where here
$S$ is the antipode of $U$. A similar treatment of $R(\chi, \rho)$
can be given based on this equation.
\section{Certain Algebras Whose Finite-Dimensional Simple
Modules Are One-Dimensional }\label{SecFinDimOne} In light of
Corollary \ref{CorSimples} we wish to consider conditions under
which the irreducible representations of an algebra are
one-dimensional with an eye towards applications of the results of
Sections \ref{SecPrelim}--\ref{SecMainHM} to certain classes of
pointed Hopf algebras. {\em Throughout this section the field} $k$
{\em is algebraically closed}. We are interested in algebras $A$
satisfying the following condition:
\begin{enumerate}
\item [] {\em $A$ is generated by an abelian group $\Gamma$ of
units of $A$ together with finitely many elements
$a_1,\dots,a_{\theta}$ and there are non-trivial characters
$\chi_1,\dots,\chi_{\theta}$ such that}
\begin{equation}\label{EqConj}
ga_i g^{-1} = \chi_i(g)a_i \text{ for all }g \in \Gamma \text{ and
}1 \leq i \leq \theta.
\end{equation}
\end{enumerate}
We denote the preceding condition by (C). Many pointed Hopf
algebras satisfy condition (C). See \cite{Irr}.

Suppose $A$ is an algebra which satisfies condition (C). We will
find sufficient conditions for all finite-dimensional simple left
$A$-modules $M$ to be one-dimensional. Finding a non-zero $m \in
M$ which satisfies $a_1{\cdot}m = \cdots = a_{\theta}{\cdot}m = 0$
is the key. The theorem of this section gives such a condition
which relates the values $\chi_j(g_i)$ to a Cartan matrix of
finite type, where $g_1, \ldots, g_{\theta} \in \Gamma$.
\begin{Lemma}\label{LemmaOneDim}
Let $A$ be an algebra satisfying (C).
\begin{enumerate}
\item[{\rm a)}] Suppose $M = km$ is a one-dimensional left
$A$-module. Then $a_i{\cdot}m = 0$ for all $1 \leq i \leq \theta$.
\end{enumerate}
Suppose that $M$ is a non-zero finite-dimensional left $A$-module.
\begin{enumerate}
\item[{\rm b)}] Assume that $M$ has a non-zero element $m$ such
that $a_i{\cdot}m = 0$ for all $1 \leq i \leq \theta$. Then $M$
contains a one-dimensional left $A$-module. \item[{\rm c)}]
Suppose that $I_1, \ldots, I_r$ partition $\{1,\dots,\theta\}$ and
$a_i$, $a_j$ skew commute whenever $i$, $j$ belong to different
$I_\ell$'s. Let $A_i$ be the subalgebra of $A$ generated by the
$a_j$'s, where  $j \in I_i$, and $\Gamma$. For all $1 \leq i \leq
\theta$ assume that non-zero finite-dimensional left $A_i$-modules
contain a one-dimensional submodule. Then $M$ contains a
one-dimensional left $A$-submodule.
\end{enumerate}
\end{Lemma}

\pf Assume the hypothesis of part {\rm a)}. Then there is a $\rho
\in \Alg(A, k)$ such that  $a{\cdot}m = \rho(a)m$ for all $a \in
A$. Let $g \in \Gamma$ and $1 \leq i \leq {\theta}$. From the
calculation
$$
\rho (ga_i)m = ga_i{\cdot}m = \chi_i(g)a_i g{\cdot}m =
\chi_i(g)\rho (a_i g)m
$$
we see that $\rho (a_i) = \chi_i (g)\rho (a_i)$. Since $\chi_i
\neq 1$ it follows that $\rho(a_i) = 0$. We have shown that
$a_i{\cdot}m = 0$ and thus part {\rm a)} is established.

As for part {\rm b)}, let $M' = \{ m \in M\,|\, a_1{\cdot}m =
\cdots = a_{\theta}{\cdot}m = 0\}$. Since $a_i g =
\chi_i(g)^{-1}ga_i$ for all $1 \leq i \leq {\theta}$ and $g \in
\Gamma$, we conclude that $M'$ is a left $\Gamma$-module. Now the
$A$-submodules of $M'$ are the $\Gamma$-submodules of the same.
Since $M'$ is finite-dimensional, $\Gamma$ is abelian, and $k$ is
algebraically closed, $M'$ contains a one-dimensional left
$\Gamma$-submodule. This concludes our proof of part {\rm b)}.

To show part {\rm c)} we may assume $r = 2$ by induction on $r$.
Thus ${\theta} > 1$ and without loss of generality we may assume
$S_1 = \{a_1, \ldots, a_s\}$ and $S_2 = \{ a_{s + 1}, \ldots,
a_{\theta}\}$ for some $1 \leq s  < {\theta}$. Since $a_i$ and
$a_j$ skew commute whenever $1 \leq i \leq s < j \leq {\theta}$,
and the elements of the commutative group $\Gamma$ skew commute
with $a_1, \ldots, a_{\theta}$, we conclude that $A_2a_i = a_i
A_2$ for all $1 \leq i \leq s$.

Let $M'$ be the set of all $m \in M$ such that $a_1{\cdot}m =
\cdots = a_s{\cdot}m = 0$. By assumption $M$ contains a
one-dimensional left $A_1$-submodule. Thus $M' \neq (0)$ by part
{\rm a)}. Since $A_2a_i = a_i A_2$ for all $1 \leq i \leq s$ it
follows that  $M'$ is a (non-zero) left $A_2$-submodule of $M$. By
assumption $M'$ contains a one-dimensional left $A_2$-submodule
$km$. Now $a_{s+1}{\cdot}m = \cdots = a_{\theta}{\cdot}m = 0$ by
part {\rm a)} again. Since $m \in M'$  by definition $a_1{\cdot}m
= \cdots = a_s{\cdot}m = 0$. Therefore $M$ contains a
one-dimensional left $A$-module by part {\rm b)}. \epf
\begin{Cor}\label{CorOneDim}
Let $A$ be an algebra satisfying condition (C). Assume that $A'$
is a subalgebra of $A$ generated by  $a_1, \ldots, a_{\theta}$ and
a subgroup $\Gamma'$ of $\Gamma$ such that the restrictions
$\chi_1|{\Gamma'}, \ldots, \chi_{\theta}|{\Gamma'} \neq 1$. Then
finite-dimensional simple left $A$-modules are one-dimensional if
the same is true for $A'$.
\end{Cor}

\pf Suppose that finite-dimensional simple left $A'$-modules are
one-dimensional and let $M$ be a finite-dimensional simple left
$A$-module. Then $M$ contains a finite-dimensional simple left
$A'$-module which must have the form $km$ by assumption. By part
{\rm a)} of Lemma \ref{LemmaOneDim} we have that $a_1{\cdot}m =
\cdots = a_{\theta}{\cdot}m = 0$. By part {\rm b)} of the same $M$
contains a one-dimensional $A$-submodule $M'$. Since $M$ is simple
$M = M'$. \epf
\medskip

Apropos of Lemma \ref{LemmaOneDim}, finite-dimensional simple left
$A$-modules are one-dimensional when $\chi_1, \ldots,
\chi_{\theta}$ are free monoid generators. For more generally:
\begin{Prop}\label{PropFinDimOneDim}
Let $A$ be an algebra satisfying condition (C). Suppose further
that $\chi_1^{k_1}\cdots \chi_{\theta}^{k_{\theta}} = 1$, where
$k_1, \ldots, k_{\theta} \geq 0$, implies $k_1 = \cdots =
k_{\theta} = 0$. Then finite-dimensional simple left $A$-modules
are one-dimensional.
\end{Prop}

\pf Let $M$ be a finite-dimensional non-zero left $A$-module.
Regarding $M$ as a left $\Gamma$-module we may write $M =
\bigoplus_{\lambda \in \widehat{\Gamma}}M_\lambda$ as the direct
sum of weight spaces, where $M_\lambda = \{ m \in M\, | \,
g{\cdot}m = \lambda (g)m \;  \forall \; g \in \Gamma\}$. Since $M$
is finite-dimensional all but finitely many of the $M_\lambda$'s
are zero. Let $1 \leq i \leq {\theta}$. Our assumption $ga_i
g^{-1} = \chi(g)a_i$ for all $g \in \Gamma$ means that
$a_i{\cdot}M_\lambda \subseteq M_{\chi_i \lambda}$ for all
$\lambda \in \widehat{\Gamma}$.

By Lemma \ref{LemmaOneDim} it suffices to show that there is a
non-zero $m \in M$ such that $a_i{\cdot}m = 0$ for all $1 \leq i
\leq {\theta}$. Suppose this is not the case. Since $M \neq (0)$
there is a $\lambda \in \widehat{\Gamma}$ such that $M_\lambda
\neq (0)$. Choose a non-zero $m \in M_{\lambda}$. By induction
there is an infinite sequence of integers $i_1, i_2, \dots$ such
that $1 \leq i_j \leq {\theta}$ for all $j \geq 1$ and
$a_{i_r}\cdots a_{i_1}{\cdot}m \neq 0$ for all $r \geq 1$. Now
$a_{i_r}\cdots a_{i_1}{\cdot}m \in M_{\chi_r\cdots
\chi_1\lambda}$. Our assumption on products of the characters
$\chi_1^{k_1}\cdots \chi_{\theta}^{k_{\theta}}$ means that
$\lambda, \chi_{i_1}\lambda, \chi_{i_2}\chi_{i_1}\lambda, \ldots$
are all distinct. But this is impossible since all but finitely
many weight spaces are zero. Therefore there is a non-zero $m \in
M$ such that $a_i{\cdot}m = 0$ for all $1 \leq i \leq {\theta}$
after all. \epf
\medskip


We recall some notions from \cite{AS2}. A {\em datum of Cartan
type}
$$
\mathcal{D} = \mathcal{D}(\Gamma, (g_i)_{1 \leq i \leq \theta},
(\chi_i)_{1 \leq i \leq \theta}, (a_{ij})_{1 \leq i,j \leq \theta})
$$
consists of an abelian group $\Gamma$, elements $g_i \in \Gamma,
\chi_i \in \widehat{\Gamma}, 1 \leq i \leq \theta,$ and a $\theta
{\times} \theta$ Cartan matrix $(a_{ij})$ satisfying
\begin{equation}\label{Cartantype}
q_{ij} q_{ji} = q_{ii}^{a_{ij}},\;q_{ii}\neq 1, \text{ with }
q_{ij} = \chi_j(g_i) \text{ for all } 1 \leq i,j \leq \theta.
\end{equation}
We define $q_i = q_{ii}$ for all $1 \leq i \leq \theta.$
Note that by \eqref{Cartantype}
\begin{equation}\label{symmetricCartan}
q_i^{a_{ij}} = q_j^{a_{ji}} \text{ for all } 1 \leq i,j \leq
\theta.
\end{equation}
Recall that a (generalized) Cartan matrix $(a_{ij})_{1 \leq i,j
\leq \theta}$ is a matrix whose entries are integers such that
$a_{ii}=2$ for all $i,$ $a_{ij} \leq 0$ for all $i \neq j$, and if
$a_{ij} = 0$ then $a_{ji}=0$ for all $i,j.$ A datum $\mathcal{D}$
of Cartan type will be called of finite Cartan type if $(a_{ij})$
is of finite type.

Let $\mathcal{D} = \mathcal{D}(\Gamma, (g_i)_{1 \leq i \leq
\theta}, (\chi_i)_{1 \leq i \leq \theta}, (a_{ij})_{1 \leq i,j
\leq \theta})$ be a datum of Cartan type. Suppose $1 \leq i,j \leq
\theta$. We say that $i$ is connected to $j$, denoted by $i \sim
j$, if there are indices $1 \leq i_1,\cdots, i_t \leq \theta,$
where  $t \geq 2,$ with $i=i_1,j=i_t$ and $a_{i_{l}i_{l+1}} \neq
0$ for all $1 \leq l <t.$ In this case we define $a(i,j) =
a_{i_1i_2} a_{i_2i_3} \cdots a_{i_{t-1}j},$ and $b(i,j) =
a_{i_2i_1} a_{i_3i_2} \cdots a_{ji_{t-1}}.$ Then it follows from
\eqref{Cartantype} that
\begin{equation}\label{aa'}
q_{i}^{a(i,j)} = q_{j}^{b(i,j)}.
\end{equation}
Connectivity is an equivalence relation and the equivalence
classes are called the {\em connected components} of $\{1,\dots,
\theta\}.$

More generally, let $R$ be a ring and $(a_{ij}) \in {\rm
M}_{\theta}(R)$ be a $\theta {\times} \theta$ matrix with
coefficients in $R$. Let $1 \leq i,j \leq \theta.$ We say that $i$
is connected to $j$ if there are indices $1 \leq i_1,\cdots, i_t
\leq \theta,$ where $t \geq 2,$ with $i=i_1,j=i_t$ and
$a_{i_{l}i_{l+1}} \neq 0$ for all $1 \leq l <t.$ In this
generality connectivity may not be an equivalence relation.

In the proof of the main theorem in this section we use the
following lemma in the special case of data of Cartan type.

\begin{Lemma}\label{connectedCartan}
Suppose $(a_{ij}) \in {\rm M}_{\theta}(\mathbb{Z})$ is a non-zero
matrix and all indices $1 \leq i,j \leq \theta$ are connected.
Assume further that $q_1,\dots,q_{\theta} \in k$ are non-zero,
that \eqref{symmetricCartan} holds, and that one of the $q_i$'s is
not a root of unity. Then:
\begin{enumerate}
\item [{\rm a)}] None of $q_1,\dots,q_{\theta}$ is a root of
unity. \item [{\rm b)}] There are roots of unity $\omega_1,
\ldots, \omega_{\theta}$ in $k,$ an element $q \in k,$ and
non-zero integers $d_1,\dots, d_{\theta}$ with $q_i = \omega_i
q^{d_i}$ for all $1 \leq i \leq \theta.$\label{oneq} \item[{\rm
c)}] $d_i a_{ij} = d_j a_{j i} \text{ for all }1 \leq i, j\leq
\theta.$\label{EqSymmetrizers} \item[{\rm d)}] Suppose $\G$ is a
group, $g_1,\dots,g_{\theta} \in \G$, and
$\chi_1,\dots,\chi_{\theta} \in \widehat{\G}$ satisfy $\chi_i(g_i)
= q_i$ and $\chi_j(g_i) \chi_i(g_j) = q_i^{a_{ij}}$ for all $1
\leq i,j \leq \theta$.
Consider the quadratic form\\
$$
Q(x_1, \ldots, x_{\theta}) = \sum_{i =1 }^{\theta} 2x_i^2d_i +
\sum_{1 \leq i < j \leq {\theta}} 2x_i x_j d_i a_{ij}.
$$
Let $k_1, \ldots , k_{\theta} \in \mathbb{Z}$ and suppose
$\chi_1^{k_1}\cdots \chi_{\theta}^{k_{\theta}} = 1$ or
$g_1^{k_1}\cdots g_{\theta}^{k_{\theta}} = 1.$ Then
$Q(k_1,\dots,k_{\theta})  =0.$
\end{enumerate}
\end{Lemma}

\pf Suppose that $q_i$ is not a root of unity and $a_{ij} \neq 0$
where $1 \leq j \leq \theta$. Since $q_i^{a_{ij}} = q_j^{a_{ji}}$,
necessarily $a_{ji} \neq 0$ and therefore $q_j$ is not a root of
unity. We now conclude that if $1 \leq \ell \leq \theta$ and $i$
is connected to $\ell$ then $q_{\ell}$ is not a root of unity. We
have shown part a) and the hypotheses of Lemma \ref{LemmaQPower}
are met. Thus part b) follows from Lemma \ref{LemmaQPower}.

By parts {\rm a)} and {\rm b)}, $q$ is not a root of unity. For $1
\leq i,j \leq \theta$ we calculate $\omega_iq^{d_ia_{ij}} =
q_i^{a_{ij}}=q_j^{a_{ji}}=\omega_jq^{d_ja_{ji}}$. Choose a
positive integer $N$ such that $\omega_i^N = 1$ for all $1 \leq i
\leq {\theta}$. Then $\displaystyle{q^{Nd_i a_{i j}} = q^{Nd_j
a_{j i}}}$. Since $q$ is not a root of unity and $N$ is a positive
integer the preceding equation implies part c).

It remains to show part d). Let $1 \leq i \leq \theta$. By
assumption $q_i^2 = \chi_i(g_i)^2 = q_i^{a_{ii}}$ and thus $a_{ii}
=2$ since $q_i$ is not a root of unity.

Now let $k_1, \ldots , k_{\theta} \in \mathbb{Z}$ and suppose that
$\chi_1^{k_1}\cdots \chi_{\theta}^{k_{\theta}} = 1$. Then
$\displaystyle{\prod_{j = 1}^{\theta}\chi_j (g_i)^{k_j} = 1}$, and
hence  $\displaystyle{\prod_{j = 1}^{\theta}\chi_j(g_i)^{k_i k_j}
= 1,}$ for all  $1 \leq i \leq {\theta}$. Suppose that
$g_1^{k_1}\cdots g_{\theta}^{k_{\theta}} = 1$. Then
$\displaystyle{\prod_{i = 1}^{\theta}\chi_j (g_i)^{k_i} = 1,}$ and
hence  $\displaystyle{\prod_{i = 1}^{\theta}\chi_j(g_i)^{k_i k_j}
= 1,}$ for all  $1 \leq j \leq {\theta}$. Thus in both cases
\begin{eqnarray*}
1
& = &  \prod_{1 \leq i,j \leq \theta} q_{ij}^{k_ik_j} \\
& = & \prod_{i = 1}^\theta q_i^{k_i^2} \prod_{1 \leq i <  j \leq
\theta}q_{ij}^{k_i k_j}
\prod_{\theta \geq i >  j \geq 1}q_{ji}^{k_i k_j} \\
& = & \prod_{i = 1}^\theta q_i^{k_i^2}
\prod_{1 \leq i <  j \leq \theta}(q_{ij}q_{ji})^{k_i k_j} \\
& = & \prod_{i = 1}^{\theta}q_i^{k_i^2} \prod_{1 \leq i <  j \leq
\theta}q_i^{a_{i  j}k_i k_j}.
\end{eqnarray*}
Raising the last expression to the $2N$ power we have $1 =
(q^N)^{Q(k_1, \ldots, k_\theta)}.$ Since $q^N$ is not a root of
unity necessarily $Q(k_1,\dots,k_{\theta}) =0.$ \epf

\medskip

We remark that part c) of the previous lemma  was shown in
\cite[Lemma 2.4]{AS3} for matrices of Cartan type in a different
way.

\begin{Lemma}\label{LemmaQPower}
Suppose that $S$ is a finite non-empty subset of non-zero elements
of $k$ which satisfies the following property: For all $x, y \in
S$ there is an $r \geq 1$, a sequence $x = x_0, x_1, \ldots, x_r =
y$ in $S$, and there are sequences $n_1, \ldots, n_r$ and
$m_0,\dots,m_{r-1}$ of non-zero integers such that
$x_{i-1}^{m_{i-1}}=x_i^{n_i}$ for all $1 \leq i \leq r$. Then
there is a $q \in k$ such that each $x \in S$ can be written as $x
= \omega q^L$ for some root of unity $\omega \in k$ and non-zero
integer $L$.
\end{Lemma}

\pf We may as well assume that there is a non-zero $x_0 \in S$.
Consider tuples ${\mathcal C} = (x_0, \ldots, x_r, n_1, \ldots,
n_r,m_0,\dots,m_{r-1})$, where the condition of the lemma is
satisfied and let $|{\mathcal C}| = n_1\cdots n_r$. By assumption
there are tuples ${\mathcal C}_1, \ldots, {\mathcal C}_s$ such
that each $x \in S$ appears as an $x_i$ in one of them.

Let $N = |{\mathcal C}_1|\cdots |{\mathcal C}_s|$. Then $N$ is a
non-zero integer. Since $k$ is algebraically closed there is a $q
\in k$ which satisfies $q^N = x_0$. Let  ${\mathcal C} = (x_0,
\ldots, x_r, n_1, \ldots, n_r,m_0,\dots,m_{r-1})$ be one of the
${\mathcal C}_i$'s. To complete the proof it suffices to show for
all $0 \leq i \leq r$ that
$$
\displaystyle{x_i = \omega'
\left(q^\frac{N}{n_1{\cdots}n_{i}}\right)^{\ell_i}}
$$
for some root of unity $\omega' \in k$ and non-zero integer
$\ell_i$. The case $i = 0$ is trivial. (By convention
$n_1{\cdots}n_i = 1$ when $i = 0$.)

Suppose that $1 < i \leq r$ and $x_{i - 1}$ has this form. Then we
can write $\displaystyle{x_{i-1} = \omega
\left(q^\frac{N}{n_1{\cdots}n_{i - 1}}\right)^{\ell_{i-1}}}$ for
some root of unity $\omega \in k$ and non-zero integer
$\ell_{i-1}$. From the calculation
$$
\displaystyle{x_i^{n_i} = x_{i - 1}^{m_{i-1}} =
\omega^{m_{i-1}}\left(q^\frac{N}{n_1\cdots
n_i}\right)^{n_{i}\ell_{i-1} m_{i-1}} =
\omega^{m_{i-1}}\left[\left(q^\frac{N}{n_1\cdots
n_i}\right)^{m_{i-1}\ell_{i-1} }\right]^{n_i}}
$$
we deduce that $x_i \neq 0$ and that $\displaystyle{x_i =
\omega'\left(q^\frac{N}{n_1\cdots
n_i}\right)^{m_{i-1}\ell_{i-1}}}$, where $\omega' \in k$ is a root
of unity. Take $l_i=m_{i-1}l_{i-1}.$ \epf
\medskip

The main result of this section is:
\begin{Theorem}\label{ThmCartan}
Let $A$ be an algebra satisfying (C) such that $\Gamma$ and
$(\chi_i)_{1 \leq i \leq \theta}$ are part of a datum $\mathcal{D}
= \mathcal{D}(\Gamma, (g_i)_{1 \leq i \leq \theta}, (\chi_i)_{1
\leq i \leq \theta}, (a_{ij})_{1 \leq i,j \leq \theta})$ of finite
Cartan type. Suppose that all $q_i,$ where $1 \leq i \leq \theta,$
are not roots of unity, and that if $1 \leq i, j \leq \theta$ are
in different connected components of $\{1,\dots,\theta\},$ then
$a_i$ and $a_j$ skew commute. Then finite-dimensional simple left
$A$-modules are one-dimensional.
\end{Theorem}

\pf Let $I_1, \ldots, I_r$ be the components of
$\{1,\dots,\theta\}$. Since the matrix $(a_{ij})_{1 \leq i,j \leq
\theta}$ is of finite type so are the matrices $(a_{ij})_{(i,j)
\in I_l {\times} I_l}$ of $I_l$ for all $1 \leq l \leq r$. Thus by
virtue of part {\rm c)} of Lemma \ref{LemmaOneDim} we may assume
$r = 1$; that is $\{1,\dots,\theta\}$ is connected.

Suppose that $\{1,\dots,\theta\}$ is connected. We will use
Proposition \ref{PropFinDimOneDim} to complete the proof. Let
$k_1, \ldots, k_{\theta} \geq 0$ and suppose that
$\chi_1^{k_1}\cdots \chi_{\theta}^{k_{\theta}} = 1$. Then $Q(k_1,
\ldots, k_{\theta}) = 0$ by part d) of Lemma
\ref{connectedCartan}. Let $b_{i j} = d_i a_{i  j}$ for all $1
\leq i, j \leq \theta$ and set $B = (b_{i  j})$. Then $B$ is a
symmetric matrix by part c) of the same lemma. Since $a_{i  i} =
2$ for all $1 \leq i \leq \theta$ we have
$$
0 = Q(k_1, \ldots, k_{\theta}) = \sum_{1 \leq i, j \leq \theta}
k_i b_{i \, j}k_j = (k_1  \cdots  k_\theta)B\left(\begin{array}{c}
k_1 \\ \vdots \\ k_\theta \end{array}\right).
$$
It will follow by Proposition \ref{PropFinDimOneDim} that all
finite-dimensional simple left $A$-modules are one-dimensional
once we show that $k_1 = \cdots = k_{\theta} = 0$.

To show the latter we follow \cite[Chapter III]{HUM}. Let
$\alpha_1,\dots,\alpha_{\theta}$ be a basis for the root system
$\Phi$ corresponding to the Cartan matrix $(a_{ij})$. Let $(,)$ be
the inner product of the euclidean vector space spanned by $\Phi$.
By definition $\displaystyle a_{ij}=
\frac{2(\alpha_i,\alpha_j)}{(\alpha_i,\alpha_i)}$. Set $x=
\sum_{i=1}^{\theta} k_i\alpha_i$. Then
$$
(x,x) = \sum_{i,j=1}^{\theta} k_ik_j(\alpha_i,\alpha_j)=\sum_{i,j=1}^{\theta}
k_ik_ja_{ij}\frac{(\alpha_i,\alpha_i)}{2}.
$$
Since $(a_{ij})$ is connected there is a non-zero $c \in
\mathbb{Q}$ such that $\displaystyle{d_i = c
\frac{(\alpha_i,\alpha_i)}{2}}$ for all $1 \leq i \leq \theta$.
Thus $0=Q(k_1,\dots,k_{\theta}) = c (x,x)$ which means $x=0$ and
consequently $k_1 = \cdots = k_{\theta} = 0$. \epf
\medskip

Let $\mathcal{D}$ be a datum of Cartan type and assume that no
$q_i$ is a root of unity. We have seen that the characters
$\chi_1,\dots,\chi_{\theta}$ are $\mathbb{Z}$-linearly independent
if $\D$ is connected and of finite type. If  $\D$ is not
connected, then linear independency fails for any non-trivial
linking. The next example shows that for connected data $\D$ the
characters are in general not linearly independent if $\D$ is not
of finite type.

\begin{Ex}\label{counter}
Let $\G$ be a free abelian group of rank two with basis $g_1,g_2,$
and $ q \in k$ not a root of unity. Let $(a_{ij}) =\begin{pmatrix}
 2& -2 \\
 -2 & 2
\end{pmatrix}.$ Define $\chi_1,\chi_2 \in \widehat{\G}$ by
$\chi_j(g_i)= q^{a_{ij}}$ for all $1 \leq i,j \leq 2.$ Then
$q_1=q_2 = q^{2}$ is not a root of unity, $\mathcal{D}(\Gamma,
(g_i)_{1 \leq i \leq 2}, (\chi_i)_{1 \leq i \leq 2}, (a_{ij})_{1
\leq i,j \leq 2})$ is a connected datum of Cartan type, but not of
finite type, and $\chi_1\chi_2=1.$
\end{Ex}
Finally we note that Theorem \ref{ThmCartan} is false if the
$q_i$'s are roots of unity.
\begin{Ex}\label{Examplesimple}
Let $N\geq 2$ be an integer, let $q$ a primitive $N$-th root of 1,
and let $A=k\langle g,a\mid g^N=1,gag^{-1} = qa\rangle$. The
algebra $A$ satisfies (C) with $\Gamma = \langle g \rangle,$ and
$\chi(g) = q,$ and $\Gamma$ and $\chi$ are part of a Cartan datum
of type $A_1.$ Let $\mathbb{Z}_N$ be the additive cyclic group of
order $N$ and $M$ be an $N$-dimensional vector space with basis
$m_i, i \in \mathbb{Z}_N.$ Then the rules $gm_i = q^im_i,$ and
$xm_i = m_{i+1}$ for all $i \in \mathbb{Z}_N$ determine an
irreducible left $A$-module structure on $M$.
\end{Ex}

\end{document}